\documentclass[a4paper]{article}
\usepackage{amsmath,amsthm,amssymb,amscd,graphicx}

\numberwithin{equation}{section}

{\theoremstyle{definition}\newtheorem{definitiona}{Definition}[section]
\newtheorem{notation}[definitiona]{Notation}

\newtheorem{remark}[definitiona]{Remark}
}
\newtheorem{proposition}[definitiona]{Proposition}
\newtheorem{lemma}[definitiona]{Lemma}

\newtheorem{corollary}[definitiona]{Corollary}

\newenvironment{definition}{\begin{definitiona}}{\mbox{} \hfill
    $\blacktriangle$ \end{definitiona}}

\newcommand{\cRtil}{\widetilde{\cR}}
\newcommand{\betatil}{\tilde{\beta}}
\newcommand{\pent}{{\text{\fontshape{n}\selectfont pent}}}
\newcommand{\coborpent}{\cobor_\pent}
\newcommand{\group}{_{\text{\fontshape{n}\selectfont group}}}
\newcommand{\coborgroup}{\cobor\group}
\newcommand{\coborhelp}{\cobor_{\text{\fontshape{n}\selectfont temp}}}
\newcommand{\hor}{_{\text{\fontshape{n}\selectfont h}}}
\newcommand{\ver}{_{\text{\fontshape{n}\selectfont v}}}
\newcommand{\lien}{\mathfrak{n}}
\newcommand{\liek}{\mathfrak{k}}
\newcommand{\liea}{\mathfrak{a}}
\newcommand{\lies}{\mathfrak{s}}
\newcommand{\Htil}{\tilde{H}}
\newcommand{\borel}{_{\mbox{\fontsize{7}{9pt}\fontshape{n}\selectfont Borel}}}
\newcommand{\lieg}{\mathfrak{g}}
\newcommand{\cont}{_{\text{\fontshape{n}\selectfont cont}}}
\newcommand{\bor}{\partial}
\newcommand{\source}{\operatorname{Source}}
\newcommand{\range}{\operatorname{Range}}
\newcommand{\Xh}{\X_{\text{\fontshape{n}\selectfont h}}}
\newcommand{\Xv}{\X_{\text{\fontshape{n}\selectfont v}}}
\newcommand{\X}{\Gamma}
\newcommand{\Y}{\Upsilon}
\newcommand{\matched}{\text{\fontshape{it}\selectfont m.p.}}
\newcommand{\funborel}{{\mathcal{F}_{\mbox{\fontsize{7}{9pt}\fontshape{n}\selectfont Borel}}}}

\newcommand{\cobor}{\operatorname{d}}
\newcommand{\coborv}{\operatorname{d}^{\text{\fontshape{n}\selectfont v}}}
\newcommand{\coborh}{\operatorname{d}^{\text{\fontshape{n}\selectfont h}}}

\newcommand{\rL}{\operatorname{L}}
\newcommand{\tetil}{\tilde{\theta}}

\newcommand{\cB}{\mathcal{B}}
\newcommand{\R}{\mathbb{R}}

\newcommand{\de}{\Delta}

\newcommand{\te}{\theta}
\newcommand{\ot}{\otimes}
\newcommand{\al}{\alpha}
\newcommand{\be}{\beta}

\newcommand{\recht}{\rightarrow}

\newcommand{\om}{\omega}
\newcommand{\B}{\operatorname{B}}
\newcommand{\si}{\sigma}
\newcommand{\io}{\iota}

\newcommand{\deh}{\hat{\Delta}}
\newcommand{\cL}{\mathcal{L}}
\newcommand{\cR}{\mathcal{R}}

\newcommand{\cW}{\mathcal{W}}

\newcommand{\cV}{\mathcal{V}}

\newcommand{\C}{\mathbb{C}}
\newcommand{\ga}{\gamma}

\newcommand{\cA}{\mathcal{A}}

\newcommand{\Z}{\mathbb{Z}}

\newcommand{\cU}{\mathcal{U}}
\newcommand{\T}{\mathbb{T}}
\newcommand{\bcU}{\overline{\mathcal{U}}}
\newcommand{\bcV}{\overline{\mathcal{V}}}
\newcommand{\bcR}{\overline{\mathcal{R}}}
\newcommand{\borh}{\partial^{\text{\fontshape{n}\selectfont h}}}
\newcommand{\borv}{\partial^{\text{\fontshape{n}\selectfont v}}}

\newcommand{\eenmatrix}[4]{\setlength{\unitlength}{1ex}
\raisebox{-1.5ex}[0ex][3.5ex]{\parbox{7.5ex}{\begin{picture}(5,5)(0,0)
\put(1.8,1.9){\line(0,1){4}}
\put(1.8,5.9){\line(1,0){4}}
\put(5.8,5.9){\line(0,-1){4}}
\put(5.8,1.9){\line(-1,0){4}}
\put(2.3,6.2){\parbox[b]{3ex}{\centering $#1$}}
\put(6.2,3.4){\parbox[b]{3ex}{$#2$}}
\put(0,3.4){\parbox[b]{1.5ex}{\raggedleft $#3$}}
\put(2.3,0){\parbox[t]{3ex}{\centering $#4$}}
\end{picture}
}}}

\newcommand{\eenmatrixmetx}[5]{\setlength{\unitlength}{1ex}
\raisebox{-1.5ex}[0ex][3.5ex]{\parbox{7.5ex}{\begin{picture}(5,5)(0,0)
\put(1.8,1.9){\line(0,1){4}}
\put(1.8,5.9){\line(1,0){4}}
\put(5.8,5.9){\line(0,-1){4}}
\put(5.8,1.9){\line(-1,0){4}}
\put(2.3,6.2){\parbox[b]{3ex}{\centering $#1$}}
\put(6.2,3.4){\parbox[b]{3ex}{$#2$}}
\put(0,3.4){\parbox[b]{1.5ex}{\raggedleft $#3$}}
\put(2.3,0){\parbox[t]{3ex}{\centering $#4$}}
\put(2.3,3.4){\parbox[b]{3ex}{\centering $#5$}}
\end{picture}
}}}

\newcommand{\kadermetx}[1]{\setlength{\unitlength}{1ex}
\raisebox{0ex}[0ex][1.5ex]{\parbox{4ex}{\begin{picture}(4,4)(0,0)
\put(0,0){\line(0,1){4}}
\put(0,4){\line(1,0){4}}
\put(4,4){\line(0,-1){4}}
\put(4,0){\line(-1,0){4}}
\put(0,1.4){\parbox[c]{4ex}{\centering $#1$}}
\end{picture}
}}}

\begin{document}
\begin{center}
{\Large\bf Measurable Kac cohomology for bicrossed products}

\bigskip

{\sc by Saad Baaj, Georges Skandalis and Stefaan
  Vaes}
\end{center}

{\footnotesize Laboratoire de Math{\'e}matiques Pures;
Universit{\'e} Blaise Pascal; B{\^a}timent de Math{\'e}matiques; F--63177 Aubi{\`e}re
Cedex (France) \\ e-mail: Saad.Baaj@math.univ-bpclermont.fr \vspace{0.3ex}\\
Institut de Math{\'e}matiques de Jussieu; Alg{\`e}bres d'Op{\'e}rateurs et
Repr{\'e}sentations; 175, rue du Chevaleret; F--75013 Paris (France) \\
e-mails: skandal@math.jussieu.fr and vaes@math.jussieu.fr}

\bigskip

\begin{abstract}\noindent
We study the Kac cohomology for matched pairs of locally compact
groups. This cohomology theory arises from the extension theory of
locally compact quantum groups.
We prove a topological version of the Kac exact
sequence and provide methods to compute the cohomology. We give
explicit calculations in several examples using results of Moore and Wigner.
\end{abstract}

\section{Introduction}
In order to construct, in a systematic way, examples of {\it finite quantum groups}, G.\ Kac developed in \cite{Kac} a method to obtain non-trivial (i.e.\
non-commutative and non-cocommutative) quantum groups as {\it extensions} of a finite group by a finite group dual. Such an extension of a finite group
$G_1$ with the dual of a finite group $G_2$ is described by the following data: \begin{itemize} \item a large group $G$ such that $G_1$ and $G_2$ are
subgroups of $G$ satisfying $G_1 \cap G_2 = \{e\}$ and $G=G_1 G_2$ (we
say that $G_1,G_2 \subset G$ is a {\it matched pair}), \item a compatible pair of
{\it $2$-cocycles} (see Equations \eqref{eq.cocycles} below). \end{itemize} Two extensions are isomorphic if and only if the the matched pairs are the same
and the pairs of $2$-cocycles are cohomologous. As such, there appears a natural cohomology group associated with a matched pair. G.~Kac found in
\cite{Kac} an {\it exact sequence} which permits to calculate this
cohomology group in terms of the usual cohomology groups of $G_1,G_2$
and $G$ with coefficients in $\T$, the group of complex numbers of
modulus $1$.

The above theory of extensions has been generalized to the framework
of locally compact quantum groups (in the sense of Kustermans and the
third author \cite{KV1,KV2}) by Vainerman and the third author \cite{VV}:
extensions in the category of locally compact quantum groups are exactly described using matched pairs with cocycles. Again, there appears a natural
cohomology group. The aim of this paper is to study this cohomology group, to prove a version of the Kac exact sequence in a locally compact setting
and to compute the cohomology in concrete examples. As such, we shall provide a precise explanation for the calculations in \cite{VV2}.

Given a matched pair of locally compact groups (see Definition \ref{def.matched}), the first two authors introduced in \cite{BS} an alternative
notion of cocycles and hence, another cohomology group. An awkward, but straightforward calculation yields that both cohomology groups agree for
matched pairs of discrete (and in particular, finite) groups. Below, we will use a more elegant approach that permits to conclude that both
cohomologies agree for the most general matched pairs of locally compact groups.

After the fundamental work of Kac \cite{Kac}, matched pairs of locally compact groups have been studied by Majid \cite{Maj1,Maj2} in order to
construct examples of Kac algebras. His definition of a matched pair of locally compact groups $G_1,G_2 \subset G$ requires that $G_1 \cap G_2 = \{e\}$
and that the multiplication map $G_1 \times G_2 \recht G$ is a homeomorphism onto $G$. The first two authors \cite{BS,BS2} gave a more general
definition, allowing $G_1G_2$ to be an open subset of $G$ with complement of measure zero. Using such a matched pair, they constructed a
multiplicative unitary, given by a pentagonal transformation. The most general definition of a matched pair has been introduced in \cite{VV} by
Vainerman and the third author and in \cite{BSV} by the authors, requiring only that $G_1 \cap G_2 = \{e\}$ and that $G_1G_2$ has a complement of
measure zero. We remark here that examples of such matched pairs, with $G_1G_2$ having empty interior, were given in \cite{BSV} and used to construct
examples of locally compact quantum groups with remarkable topological properties.

We mention that algebraic work on matched pairs and Kac cohomology for
Hopf algebras and Lie algebras has been done in e.g.\ \cite{andrus,Mas,schauenburg}.

\section{Preliminaries}
In this paper, all locally compact spaces will be supposed {\it second
  countable}. We denote by $\T$ the group of complex numbers of
  modulus $1$, which we will often write additively through the
  identification with $\R/\Z$.

\begin{definition} \label{def.matched}
We call $G_1,G_2 \subset G$ a {\it matched pair} of locally compact
groups if $G_1,G_2$ are closed subgroups of the locally compact group
$G$ such that $G_1 \cap G_2 =\{e\}$ and $G \setminus G_1G_2$ has Haar
measure zero.
\end{definition}

\begin{notation}
Given a matched pair $G_1,G_2 \subset G$ of locally compact groups, we
define, almost everywhere on $G$,
\begin{align*}
p_1 : G \recht G_1 \quad, \quad p_2 : G \recht G_2 \quad &\text{such
  that}\quad x = p_1(x) \, p_2(x) \; ,\\
q_1 : G \recht G_1 \quad, \quad q_2 : G \recht G_2 \quad &\text{such
  that}\quad x = q_2(x) \, q_1(x) \; .
\end{align*}
\end{notation}

We will deal with cohomology theories with coefficients in Polish
$G$-modules. Therefore, the following will be useful to us.

\begin{notation} \label{not.Lspace}
Let $X$ be a standard Borel space equiped with a Borel measure class. Let $A$ be a Polish space. We define $\rL(X,A)$ to be the set of equivalence
classes of Borel measurable functions of $X$ to $A$ identifying functions equal almost everywhere. Choosing a finite measure $\mu$ on $X$ in the
given measure class and a bounded complete metric $\rho$ on $A$, we can define
$$\rho_\mu(F,G) := \int \rho(F(x),G(x)) \; d\mu(x) \quad\text{for all}\;\; F,G \in \rL(X,A) \; .$$
In this way, $\rL(X,A)$ is a Polish space. The topology on $\rL(X,A)$
does not depend on the choice of $\mu$ or $\rho$, see the Corollary to
Proposition 6 in \cite{CCM}.
\end{notation}

In Theorem 1 of \cite{CCM}, the Fubini theorem is used to prove that
there are natural isomorphisms $\rL(X,\rL(Y,A)) \cong \rL(X \times
Y,A) \cong \rL(Y,\rL(X,A))$ for all Polish spaces $A$.

We finally introduce the {\it measurable cohomology} of a locally
compact group $G$, as studied by Moore \cite{CCM} and D.\ Wigner
\cite{DW}.

Let $G$ be a locally compact group and $A$ a Polish $G$-module. We
write face operators
$$\bor_i : G^{n+1} \recht G^{n} : \bor_i(g_0,\ldots,g_n) = \begin{cases}
(g_1,\ldots,g_n) \quad\text{if}\;\; i=0, \\
(g_0,\ldots,g_{i-1} g_{i},\ldots,g_n) \quad\text{if}\;\; i=1,\ldots,n
,\\
(g_0,\ldots,g_{n-1}) \quad\text{if}\;\; i=n+1 \; .
\end{cases}$$
Dualizing, we can write at least two natural measurable cochain
complexes. First, consider
\begin{align}
& \cobor_i : \rL(G^n,A) \recht \rL(G^{n+1},A) : (\cobor_i
F)(\overrightarrow{g}) = \begin{cases} g_0 \cdot F(\bor_0
  \overrightarrow{g}) \quad\text{if}\;\; i=0 \; , \notag \\
F(\bor_i \overrightarrow{g}) \quad\text{if}\;\; i=1,\ldots n+1 \; ,
\end{cases} \text{where}\; \overrightarrow{g}=(g_0,\ldots,g_n)\; ,\\
& \cobor : \rL(G^n,A) \recht \rL(G^{n+1},A) : \cobor =
\sum_{i=0}^{n+1} (-1)^i \cobor_i \; .\label{eq.definitioncoboundary}
\end{align}

\begin{definition}
The cohomology of the cochain complex $(\rL(G^n,A))_n$ defined above
is denoted by $H(G,A)$ and called the {\it measurable cohomology} of
the locally compact group $G$ with coefficients in the Polish
$G$-module $A$.
\end{definition}

Instead of using $\rL(G^n,A)$, we can use the $\Z$-module
$\funborel(G^n,A)$ of Borel measurable functions from $G^n$ to $A$. We
define the coboundary with the same formula as in Equation
\eqref{eq.definitioncoboundary}. Moore \cite{CCM} proved that the
obvious cochain transformation $(\funborel(G^n,A)) \recht
(\rL(G^n,A))$ is a cohomology isomorphism.

In fact, there is more. Let $G$ be a locally compact group and
consider a certain category of $G$-modules. Suppose that $H(G,A)$ and
$\Htil(G,A)$ are two cohomology theories satisfying the {\it Buchsbaum
  criterion} \cite{Buchs}:
\begin{itemize}
\item every short exact sequence of $G$-modules, gives rise, in a
  natural way, to a long exact cohomology sequence,
\item effaceability, i.e.\ for every $\alpha \in H^n(G,A)$ there exists a short exact
  sequence $0 \recht A \recht B \recht C \recht 0$ such that $\alpha$
  is $0$ in $H^n(G,B)$.
\end{itemize}
If now $H^0(G,A)$ and $\Htil^0(G,A)$ are naturally isomorphic, we can
conclude that $H^n(G,A)$ and $\Htil^n(G,A)$ are naturally isomorphic
for all $n$. Further, any natural sequence of homomorphisms $H^n(G,A)
\recht \Htil^n(G,A)$ which is {\it connected} (i.e.\ respects the long
exact cohomology sequences) and which gives an isomorphism for $n=0$,
will be an isomorphism for all $n$.

Another way of describing the measurable cohomology of $G$ is through
the right notion of a \emph{free resolution}.

\begin{definition}
We say that a Polish $G$-module $A$ is \emph{free} if there exists a
Polish $G$-module $B$ such that $A \cong \rL(G,B)$.

Let $A$ be a Polish $G$-module. We call
\begin{equation} \label{eq.resolution}
0 \longrightarrow A \longrightarrow A_0 \longrightarrow A_1
\longrightarrow \cdots
\end{equation}
a resolution of $A$, if all $A_i$ are Polish $G$-modules, the arrows
are $G$-equivariant and continuous and the sequence in the previous
equation is exact.

We say that the resolution in Equation \eqref{eq.resolution} is a
\emph{free resolution} of $A$ if every Polish $G$-module $A_i$, $i
\geq 0$ is free.
\end{definition}

Whenever $0 \longrightarrow A \longrightarrow A_0 \longrightarrow A_1
\longrightarrow \cdots$ is a free resolution of the Polish $G$-module
$A$, the measurable cohomology $H(G,A)$ is the cohomology of the
complex
$$A_0^G \longrightarrow A_1^G \longrightarrow \cdots$$
where $A_i^G$ denotes the $G$-fixed points of $A_i$. Again, as above, if we have two free resolutions of $A$:
$$0 \longrightarrow A \longrightarrow A_0 \longrightarrow A_1
\longrightarrow \cdots \quad\text{and}\quad 0 \longrightarrow A \longrightarrow B_0 \longrightarrow B_1
\longrightarrow \cdots$$
and if we are given continuous $G$-equivariant homomorphisms $A_i \recht B_i$ intertwining the two free resolutions, then these homomorphisms induce a cohomology isomorphism.

In Section \ref{sec.computational}, we explain the methods developed
by D.\ Wigner to compute the measurable cohomology $H(G,A)$ in certain cases.

\section{Two Kac $\mathbf{2}$-cohomology groups and the Kac bicomplex}

The Kac $2$-cohomology appears in two natural ways. In both pictures,
it is at first somehow awkward to write the cocycle relations. We will
see below, how a much more natural way of writing these relations can
be obtained. This will allow also to unify both pictures and prove
that they define the same $2$-cohomology group.

\subsection{$\mathbf{2}$-cohomology of matched pairs} \label{subsec.bicros}

First of all, we
want to classify extensions
\begin{equation}\label{eq.extension}
e \longrightarrow (\rL^\infty(G_1),\de_1) \longrightarrow (M,\de)
\longrightarrow (\cL(G_2),\deh_2) \longrightarrow e
\end{equation}
where $(M,\de)$ is a locally compact quantum group. Here, we do not explain the
notion of such an extension in the framework of locally compact quantum groups
(see Definition 3.2 in \cite{VV}). All good (i.e.,
{\it cleft}) extensions can be written as a cocycle bicrossed product
of $G_1$ and $G_2$, see Theorem 3.6 in \cite{VV}. Using Remark 5.3 in
\cite{BSV} and Lemma 4.11 in \cite{VV}, this means that
\begin{itemize}
\item there exists a locally compact group $G$ such that
$G_1,G_2$ are closed subgroups of $G$ forming a matched pair in the
sense of Definition \ref{def.matched},
\item there exists a compatible pair $(\cU,\cV)$ of $2$-cocycles on
  the matched pair $G_1,G_2$. This means that $\cU : G_2 \times G_1
  \times G_1 \recht \T$ and $\cV : G_2 \times G_2 \times G_1 \recht
  \T$ are measurable maps satisfying the (awkward) relations
\begin{align}
& \cU(p_2(sg),h,k) \; \bcU(s,gh,k) \; \cU(s,g,hk) \; \bcU(s,g,h) = 1
\; , \notag \\
& \cV(t,r,g) \; \bcV(st,r,g) \; \cV(s,tr,g) \; \bcV(s,t,p_1(rg)) = 1
\; , \label{eq.cocycles} \\
& \cU(t,g,h) \; \bcU(st,g,h) \; \cU(s,p_1(tg),p_1(p_2(tg)h)) \;
\cV(p_2(sp_1(tg)),p_2(tg),h) \; \bcV(s,t,gh)  \; \cV(s,t,g) = 1 \; , \notag
\end{align}
for almost all $s,t,r \in G_2$ and $g,h,k \in G_1$.
\end{itemize}
The locally compact quantum group $(M,\de)$ is the cocycle bicrossed
product constructed with this data of a matched pair and a compatible
pair of $2$-cocycles.

From Proposition 3.8 in \cite{VV}, we know that
two extensions are isomorphic if and only if the associated matched
pairs of locally compact groups are the same and there exists a measurable map
$\cR : G_2 \times G_1 \recht \T$ such that the pairs of $2$-cocycles
$(\cU,\cV)$ differ by a trivial pair of $2$-cocycles
$(\cU_\cR,\cV_\cR)$ defined by
\begin{equation}\label{eq.trivialcocycles}
\cU_\cR(s,g,h) = \bcR(p_2(sg),h) \; \cR(s,gh) \; \bcR(s,g) \quad , \quad
\cV_\cR(s,t,g) = \cR(t,g) \; \bcR(st,g) \; \cR(s,p_1(tg)) \; .
\end{equation}
This leads us to the following definition (Terminology 4.21 in \cite{VV}).
\begin{definition}\label{def.Kaccohomology}
Let $G_1,G_2 \subset G$ be a matched pair of locally compact
groups. The associated {\it group of extensions} is defined as the group of
pairs of $2$-cocycles $(\cU,\cV) \in \rL(G_2 \times G_1^2,\T) \oplus
\rL(G_2^2 \times G_1,\T)$
satisfying Equation \eqref{eq.cocycles}, divided by the subgroup of trivial cocycles defined by Equation \eqref{eq.trivialcocycles}.
\end{definition}
As a conclusion, we see that extensions \eqref{eq.extension} are
classified by a matched pair $G_1,G_2 \subset G$ and the associated
group of extensions.

\subsection{$\mathbf{2}$-cohomology of pentagonal transformations}

Secondly, fix a matched pair $G_1,G_2 \subset G$. There is
an associated bicrossed product locally compact quantum group $(M,\de)$ with
multiplicative unitary $W$. By definition (see \cite{BS}), a multiplicative
unitary is a unitary operator on a tensor square $H \ot H$ of a
Hilbert space $H$, satisfying the {\it pentagonal equation}
$$W_{12} \, W_{13} \, W_{23} = W_{23} \, W_{12}$$
on $H \ot H \ot H$. For the bicrossed product of $G_1,G_2 \subset G$,
this multiplicative unitary is given, as in \cite{BS,BS2} and up
to some identifications, by the formula $(W \xi)(x,y) = d(x,y)^{1/2}
\xi(w(x,y))$ for $x,y \in G$, where
\begin{equation} \label{eq.pentagonaltransformation}
w(x,y) = (x
p_1(p_2(x)^{-1}y),p_2(x)^{-1}y)
\end{equation}
and where $d(x,y)$ is the
Radon-Nikodym derivative making $W$ a unitary operator on $L^2(G
\times G)$. In D{\'e}finition 8.24 of \cite{BS}, a $2$-cocycle for the
matched pair $G_1,G_2 \subset G$ is defined as a measurable function
$\tetil : G \times G \recht \T$ such that $W_{\tetil} := \tetil \, W$ is still a
multiplicative unitary. Here, $\tetil$ denotes as well the multiplication
operator by the function $\tetil$. A trivial $2$-cocycle is a $2$-cocycle
of the form $\tetil(x,y) = t(x,y) \overline{t(w(x,y))}$, where $t(x,y) =
a(x) \, a(y)$ for some measurable function $a : G \recht \T$. In that
case $W_{\tetil} = (a \ot a) W (a^* \ot a^*)$, which motivates why such a
$\tetil$ is considered to be trivial. Dividing the group of $2$-cocycles by trivial
$2$-cocycles, we get again a $2$-cohomology group. For finite groups
$G_1,G_2$, the computational
argument in Section 4.4 of \cite{VV} allows to conclude that this
$2$-cohomology group is isomorphic with the $2$-cohomology group
defined above using pairs $(\cU,\cV)$.
The more delicate general case
will be dealt with below.

Observing that $w = w_1 \circ w_2$, where
$w_1(x,y) = (x p_1(y) , y)$ and $w_2(x,y) = (x,p_2(x)^{-1}y)$, we
write $W = W_2 W_1$. To simplify formulas, we may as well define a
$2$-cocycle as a measurable function $\te : G \times G \recht \T$
such that $W_2 \, \te \, W_1$ is a multiplicative unitary. Of course, one
can pass from $\te$ to $\tetil$ by the formula $\tetil(x,y) =
\te(w_2(x,y))$. The $2$-cocycle relation for $\te$ becomes
\begin{equation}\label{eq.pentagonalcocycles}
\te(x,y) \; \te(x p_1(y), p_2(y) z) \; \te(y,z) = \te(p_2(x) y, z)
\; \te(x,y p_1(z))
\end{equation}
for almost all $x,y,z \in G$. Trivial $2$-cocycles are given by the
formula
\begin{equation}\label{eq.trivialpentagonalcocycles}
\te(x,y) = a(x) \; a(p_2(x)y) \; \overline{a(xp_1(y))} \;
\overline{a(y)}
\end{equation}
for some measurable function $a : G \recht \T$ and almost all $x,y \in G$.

\begin{definition} \label{def.pentagonalcohomology}
Let $G_1,G_2 \subset G$ be a matched pair. The {\it $2$-cohomology group
associated with the pentagonal transformation}
\eqref{eq.pentagonaltransformation} is defined as the group of
cocycles $\te \in \rL(G \times G,\T)$ satisfying Equation
\eqref{eq.pentagonalcocycles}, divided by the subgroup of trivial
cocycles defined by Equation \eqref{eq.trivialpentagonalcocycles}.
\end{definition}

\subsection{The Kac bicomplex}

Fix now a matched pair $G_1,G_2 \subset G$. Define the closed subspace
$\X_{11} \subset G_1 \times G_1 \times G_2 \times G_2$ as follows:
$$\X_{11}:= \Bigl\{ \eenmatrix{s}{g}{h}{t} \Big| \; g,h \in G_1, s,t
\in G_2, \; sg = ht \Bigr\} \; .$$

\begin{lemma} \label{lemma.measureclass}
The maps
$$\X_{11} \recht G_1 \times G_2 : \eenmatrix{s}{g}{h}{t} \mapsto
\begin{cases} (g,s) \\ (g,t) \\ (h,s) \\ (h,t) \end{cases}
\qquad\text{and}\qquad
\X_{11} \recht G : \eenmatrix{s}{g}{h}{t} \mapsto sg$$
are injective. Their ranges have complement of measure zero and all
these maps define the same measure class on $\X_{11}$.
\end{lemma}
\begin{proof}
This follows immediately from Proposition 3.2 in \cite{BSV} and the
remarks following that proposition.
\end{proof}

When $x \in G_1G_2 \cap G_2G_1$ (and, as follows from the previous lemma,
almost all $x \in G$ are like
that), we sometimes
write
$\raisebox{0ex}[3.5ex][0ex]{\eenmatrixmetx{s}{g}{h}{t}{x}}$ to denote the element
$\eenmatrix{s}{g}{h}{t} \in \X_{11}$ satisfying $sg = ht = x$. We even
use $\kadermetx{x}$ to denote the same element.

So, we defined the space $\X_{11}$ by labelling the edges of a
square. We will use this to give a {\it non-equivariant} image of our
cohomology theory. There is an analogous {\it equivariant} image. We
define $\Y_{11} \subset G^4$ as follows:
$$\Y_{11}:= \Bigl\{
\raisebox{0.5ex}[3.5ex][3ex]{\parbox{8ex}{\setlength{\unitlength}{1ex}\begin{picture}(7,7)(-0.5,0)
\put(1,1){\line(1,0){4}}
\put(5,1){\line(0,1){4}}
\put(1,5){\line(1,0){4}}
\put(1,1){\line(0,1){4}}
\put(-0.3,-0.3){\mbox{$z$}}
\put(-0.3,5.3){\mbox{$x$}}
\put(5.3,5.3){\mbox{$y$}}
\put(5.3,-0.3){\mbox{$w$}}
\end{picture}}} \Big| \; x,y,z,w \in G, \;\; x^{-1}z, y^{-1}w \in G_1,
\;\;  x^{-1}y, z^{-1}w \in G_2 \Bigr\} \; .$$
There is a natural action of $G$ on $\Y_{11}$ given by $$v \cdot \Biggl(
\raisebox{0.5ex}[0ex][3.5ex]{\parbox{8ex}{\setlength{\unitlength}{1ex}\begin{picture}(7,7)(-0.5,0)
\put(1,1){\line(1,0){4}}
\put(5,1){\line(0,1){4}}
\put(1,5){\line(1,0){4}}
\put(1,1){\line(0,1){4}}
\put(-0.3,-0.3){\mbox{$z$}}
\put(-0.3,5.3){\mbox{$x$}}
\put(5.3,5.3){\mbox{$y$}}
\put(5.3,-0.3){\mbox{$w$}}
\end{picture}}} \Biggr) =
\raisebox{0.5ex}[0ex][3.5ex]{\parbox{8ex}{\setlength{\unitlength}{1ex}\begin{picture}(7,7)(-0.5,0)
\put(1,1){\line(1,0){4}}
\put(5,1){\line(0,1){4}}
\put(1,5){\line(1,0){4}}
\put(1,1){\line(0,1){4}}
\put(-0.7,-0.3){\mbox{$vz$}}
\put(-0.7,5.3){\mbox{$vx$}}
\put(4.5,5.3){\mbox{$vy$}}
\put(4.5,-0.3){\mbox{$vw$}}
\end{picture}}}$$
and a natural homeomorphism $G \times \X_{11} \recht \Y_{11}$ given by
$$\Biggl( \; v \; , \; \; \eenmatrixmetx{s}{g}{h}{t}{x} \Biggr) \mapsto
\raisebox{0.5ex}[0ex][3.5ex]{\parbox{8ex}{\setlength{\unitlength}{1ex}\begin{picture}(7,7)(-0.5,0)
\put(1,1){\line(1,0){4}}
\put(5,1){\line(0,1){4}}
\put(1,5){\line(1,0){4}}
\put(1,1){\line(0,1){4}}
\put(-0.7,-1){\mbox{$vh$}}
\put(-0.3,5.3){\mbox{$v$}}
\put(4.5,5.3){\mbox{$vs$}}
\put(4.5,-1){\mbox{$vx$}}
\end{picture}}}$$

We define more generally the space $\X_{pq}$ (which is, stricto sensu,
a closed subspace of $G_1^{p(q+1)} \times G_2^{(p+1)q}$) consisting of
elements $X \in \X_{pq}$ defined by
\begin{equation}\label{eq.generalelement}
\setlength{\unitlength}{1.5ex} X =
\parbox[c]{34.5ex}{\begin{picture}(23,23)(0,0)
\put(2,2){\line(1,0){10}}
\put(13,1.6){\parbox{4.5ex}{\centering $\ldots$}}
\put(17,2){\line(1,0){6}}
\put(2,6){\line(1,0){10}}
\put(13,5.6){\parbox{4.5ex}{\centering $\ldots$}}
\put(17,6){\line(1,0){6}}
\put(2,14){\line(1,0){10}}
\put(13,13.6){\parbox{4.5ex}{\centering $\ldots$}}
\put(17,14){\line(1,0){6}}
\put(2,18){\line(1,0){10}}
\put(13,17.6){\parbox{4.5ex}{\centering $\ldots$}}
\put(17,18){\line(1,0){6}}
\put(2,22){\line(1,0){10}}
\put(13,21.6){\parbox{4.5ex}{\centering $\ldots$}}
\put(17,22){\line(1,0){6}}
\put(2,2){\line(0,1){6}}
\put(1.8,9.3){\mbox{$\vdots$}}
\put(2,12){\line(0,1){10}}
\put(6,2){\line(0,1){6}}
\put(5.8,9.3){\mbox{$\vdots$}}
\put(6,12){\line(0,1){10}}
\put(10,2){\line(0,1){6}}
\put(9.8,9.3){\mbox{$\vdots$}}
\put(10,12){\line(0,1){10}}
\put(19,2){\line(0,1){6}}
\put(18.8,9.3){\mbox{$\vdots$}}
\put(19,12){\line(0,1){10}}
\put(23,2){\line(0,1){6}}
\put(22.8,9.3){\mbox{$\vdots$}}
\put(23,12){\line(0,1){10}}
\put(2,22.5){\parbox[b]{6ex}{\centering $s_{01}$}}
\put(6,22.5){\parbox[b]{6ex}{\centering $s_{02}$}}
\put(19,22.5){\parbox[b]{6ex}{\centering $s_{0q}$}}
\put(2,18.5){\parbox[b]{6ex}{\centering $s_{11}$}}
\put(6,18.5){\parbox[b]{6ex}{\centering $s_{12}$}}
\put(19,18.5){\parbox[b]{6ex}{\centering $s_{1q}$}}
\put(2,14.5){\parbox[b]{6ex}{\centering $s_{21}$}}
\put(6,14.5){\parbox[b]{6ex}{\centering $s_{22}$}}
\put(19,14.5){\parbox[b]{6ex}{\centering $s_{2q}$}}
\put(2,0.7){\parbox[b]{6ex}{\centering $s_{p1}$}}
\put(6,0.7){\parbox[b]{6ex}{\centering $s_{p2}$}}
\put(19,0.7){\parbox[b]{6ex}{\centering $s_{pq}$}}
\put(-0.1,3.7){\mbox{$g_{p0}$}}
\put(-0.1,15.7){\mbox{$g_{20}$}}
\put(-0.1,19.7){\mbox{$g_{10}$}}
\put(5,3.7){\mbox{$g_{p1}$}}
\put(5,15.7){\mbox{$g_{21}$}}
\put(5,19.7){\mbox{$g_{11}$}}
\put(10.2,3.7){\mbox{$g_{p2}$}}
\put(10.2,15.7){\mbox{$g_{22}$}}
\put(10.2,19.7){\mbox{$g_{12}$}}
\put(23.2,3.7){\mbox{$g_{pq}$}}
\put(23.2,15.7){\mbox{$g_{2q}$}}
\put(23.2,19.7){\mbox{$g_{1q}$}}
\end{picture}}
\end{equation}
where all $g_{ij} \in G_1$, $s_{ij} \in G_2$ and every small square of
the above picture belongs to $\X_{11}$. This means, e.g., that
$s_{12}g_{22} = g_{21}s_{22}$. More generally, this means that, if you
choose two vertices in $X$ and a path between them, then the result of
the multiplication of all the edges along the path, does not depend on the
chosen path. As such, we define $\X_{pq}$ whenever $p+q \neq 0$. We
define $\X_{00}$ to be one point and remark that $\X_{p0}=G_1^p$ and $\X_{0q}=G_2^q$.

Of course, we have again an analogous equivariant image $\Y_{pq}
\subset G^{(p+1)\times (q+1)}$, consisting of $(p+1) \times (q+1)$-matrices $Y$
with entries in $G$
\begin{equation}\label{eq.generalelementY}
\setlength{\unitlength}{1.5ex}
Y = \parbox[c]{34.5ex}{\begin{picture}(23,23)(0,0)
\put(2,2){\line(1,0){10}}
\put(13,1.6){\parbox{4.5ex}{\centering $\ldots$}}
\put(17,2){\line(1,0){6}}
\put(2,6){\line(1,0){10}}
\put(13,5.6){\parbox{4.5ex}{\centering $\ldots$}}
\put(17,6){\line(1,0){6}}
\put(2,14){\line(1,0){10}}
\put(13,13.6){\parbox{4.5ex}{\centering $\ldots$}}
\put(17,14){\line(1,0){6}}
\put(2,18){\line(1,0){10}}
\put(13,17.6){\parbox{4.5ex}{\centering $\ldots$}}
\put(17,18){\line(1,0){6}}
\put(2,22){\line(1,0){10}}
\put(13,21.6){\parbox{4.5ex}{\centering $\ldots$}}
\put(17,22){\line(1,0){6}}
\put(2,2){\line(0,1){6}}
\put(1.8,9.3){\mbox{$\vdots$}}
\put(2,12){\line(0,1){10}}
\put(6,2){\line(0,1){6}}
\put(5.8,9.3){\mbox{$\vdots$}}
\put(6,12){\line(0,1){10}}
\put(10,2){\line(0,1){6}}
\put(9.8,9.3){\mbox{$\vdots$}}
\put(10,12){\line(0,1){10}}
\put(19,2){\line(0,1){6}}
\put(18.8,9.3){\mbox{$\vdots$}}
\put(19,12){\line(0,1){10}}
\put(23,2){\line(0,1){6}}
\put(22.8,9.3){\mbox{$\vdots$}}
\put(23,12){\line(0,1){10}}
\put(0.8,22.5){\mbox{$x_{00}$}}
\put(5,22.5){\mbox{$x_{01}$}}
\put(9,22.5){\mbox{$x_{02}$}}
\put(22,22.5){\mbox{$x_{0,q}$}}
\put(-0.3,18.3){\mbox{$x_{10}$}}
\put(23.3,18.3){\mbox{$x_{1q}$}}
\put(3.8,18.3){\mbox{$x_{11}$}}
\put(7.8,18.3){\mbox{$x_{12}$}}
\put(-0.3,14.3){\mbox{$x_{20}$}}
\put(3.8,14.3){\mbox{$x_{21}$}}
\put(7.8,14.3){\mbox{$x_{22}$}}
\put(23.3,14.3){\mbox{$x_{2q}$}}
\put(0.8,1){\mbox{$x_{p0}$}}
\put(5,1){\mbox{$x_{p1}$}}
\put(9,1){\mbox{$x_{p2}$}}
\put(22,1){\mbox{$x_{p,q}$}}
\end{picture}}
\end{equation}
such that the elements on a fixed row define the
same element of $G/G_2$ and the elements on a fixed column define the same element of $G/G_1$. More formally, $x_{ij}^{-1}
x_{kj} \in G_1$ and $x_{ij}^{-1} x_{ik} \in G_2$.

Again, we have an action of $v \in G$ on $Y \in \Y_{pq}$, multiplying
all $x_{ij}$ in Equation \eqref{eq.generalelementY} by $v$ on the
left. We get a homeomorphism $\Y_{pq} \recht G \times \X_{pq}$, which
sends an element $Y$ to the couple $(x_{00},X)$, where $X$ is defined
by $g_{ij} = x_{i-1,j}^{-1} x_{ij}$ and $s_{ij} = x_{i,j-1}^{-1}x_{ij}$.

\begin{remark} \label{rem.measurespace}
An element of $\X_{pq}$ is uniquely determined once you know one row and one column. The mapping $\X_{pq} \recht G_1^p \times G_2^q$ picking out one column and one row is injective
and its image has a complement of measure zero.
It follows from Lemma \ref{lemma.measureclass} that
all these mappings induce the same measure class on $\X_{pq}$.

Also the mapping $\X_{pq} \recht \X_{k,l} \times \X_{p-k,q-l}$, sending a matrix to its upper left and lower right corner, is injective. Its image has a complement of measure zero
and the map is a measure class isomorphism.
\end{remark}

We define horizontal and vertical face operators on $\Y_{pq}$:
\begin{align*} \borh_i : \Y_{pq} \recht \Y_{p,q-1} : \borh_i \;\;
  &\text{removes the} \; i\text{-th column}\quad(\text{for}\;\; i = 0,
\ldots , q) \quad\text{and} \\ \borv_j : \Y_{pq} \recht \Y_{p-1,q} :
\borv_j \;\;
  &\text{removes the} \; j\text{-th row}
\quad(\text{for}\;\; j = 0, \ldots , p) \; .
\end{align*}

These face operators are obviously $G$-equivariant and so, we get
face operators on $\X_{pq}$, defined as follows.
The face $\borh_i$ contracts the $i$-th column multiplying the
adjacent horizontal edges. In a concrete example, this means that
\setlength{\unitlength}{1ex}
$$\borh_0 \Bigl( \parbox[c]{12ex}{\begin{picture}(12,8)(0,0)
\put(2,2){\line(1,0){8}}
\put(2,6){\line(1,0){8}}
\put(2,2){\line(0,1){4}}
\put(6,2){\line(0,1){4}}
\put(10,2){\line(0,1){4}}
\put(0.5,3.4){\mbox{$g$}}
\put(4.5,3.4){\mbox{$h$}}
\put(10.2,3.4){\mbox{$k$}}
\put(2,0){\parbox[b]{4ex}{\centering $t$}}
\put(6,0){\parbox[b]{4ex}{\centering $t'$}}
\put(2,6.4){\parbox[b]{4ex}{\centering $s$}}
\put(6,6.4){\parbox[b]{4ex}{\centering $s'$}}
\end{picture}}\Bigr) = \eenmatrix{s'}{k}{h}{t'}\quad , \quad
\borh_1 \Bigl( \parbox[c]{12ex}{\begin{picture}(12,8)(0,0)
\put(2,2){\line(1,0){8}}
\put(2,6){\line(1,0){8}}
\put(2,2){\line(0,1){4}}
\put(6,2){\line(0,1){4}}
\put(10,2){\line(0,1){4}}
\put(0.5,3.4){\mbox{$g$}}
\put(4.5,3.4){\mbox{$h$}}
\put(10.2,3.4){\mbox{$k$}}
\put(2,0){\parbox[b]{4ex}{\centering $t$}}
\put(6,0){\parbox[b]{4ex}{\centering $t'$}}
\put(2,6.4){\parbox[b]{4ex}{\centering $s$}}
\put(6,6.4){\parbox[b]{4ex}{\centering $s'$}}
\end{picture}}\Bigr)
= \eenmatrix{ss'}{k}{g}{tt'}
\quad\quad\text{and}\quad\quad
\borh_2 \Bigl( \parbox[c]{12ex}{\begin{picture}(12,8)(0,0)
\put(2,2){\line(1,0){8}}
\put(2,6){\line(1,0){8}}
\put(2,2){\line(0,1){4}}
\put(6,2){\line(0,1){4}}
\put(10,2){\line(0,1){4}}
\put(0.5,3.4){\mbox{$g$}}
\put(4.5,3.4){\mbox{$h$}}
\put(10.2,3.4){\mbox{$k$}}
\put(2,0){\parbox[b]{4ex}{\centering $t$}}
\put(6,0){\parbox[b]{4ex}{\centering $t'$}}
\put(2,6.4){\parbox[b]{4ex}{\centering $s$}}
\put(6,6.4){\parbox[b]{4ex}{\centering $s'$}}
\end{picture}}\Bigr)
= \eenmatrix{s}{h}{g}{t} \;\; .
$$
It will be clear how to define the vertical face $\borv_j$ contracting
the $j$-th row and multiplying the adjacent vertical edges. Observe
that $\borh_i \borh_j =\borh_{j-1} \borh_i$ if $i < j$.

\bigskip

Let $A$ be a Polish $G$-module.
\begin{notation} \label{not.identification}
Following Notation \ref{not.Lspace}, we consider the Polish $\Z$-modules
$\rL(\X_{pq},A)$ and $\rL(\Y_{pq},A)$. Using the action of $G$ on
$\Y_{pq}$, we turn $\rL(\Y_{pq},A)$ into a Polish $G$-module,
defining $$(x \cdot F)(Y) = x \cdot F(x^{-1} \cdot Y)$$ for $x \in G$
and $Y \in \Y_{pq}$.
We denote by $\rL(\Y_{pq},A)^G$ the $\Z$-module of $G$-invariant elements
of $\rL(\Y_{pq},A)$. We have a natural identification $\rL(\X_{pq},A)
\cong \rL(\Y_{pq},A)^G$.
\end{notation}

We define a bicomplex (see \cite{maclane})
\begin{equation}\label{eq.Kacbicomplex}
\begin{CD}
\vdots @. \vdots @. \\
@A{\coborv}AA @A{\coborv}AA @. \\
\rL(\X_{21},A) @>{\coborh}>> \rL(\X_{22},A) @>{\coborh}>> \cdots\\
@A{\coborv}AA @A{\coborv}AA @. \\
\rL(\X_{11},A) @>{\coborh}>> \rL(\X_{12},A) @>{\coborh}>> \cdots
\end{CD}
\end{equation}

where the arrows can be defined most easily using Notation
\ref{not.identification} and the equivariant coboundary operators
\begin{equation} \label{eq.equivariant}
\begin{split}
\coborh : \rL(\Y_{pq},A) \recht \rL(\Y_{p,q+1},A) : (\coborh F)(Y) &=
\sum_{i = 0}^{q+1} (-1)^i \; F(\borh_i Y) \; , \\
\coborv : \rL(\Y_{pq},A) \recht \rL(\Y_{p+1,q},A) : (\coborv F)(Y) &=
\sum_{j = 0}^{p+1} (-1)^j \; F(\borv_j Y) \; .
\end{split}
\end{equation}
A small calculation reveals that, on $\X_{pq}$ rather than $\Y_{pq}$,
we get
\begin{align*}
\coborh : \rL(\X_{pq},A) \recht \rL(\X_{p,q+1},A) : (\coborh F)(X) &=
s_{01} \cdot F(\borh_0 X) + \sum_{i = 1}^{q+1} (-1)^i \; F(\borh_i X) \;
, \\
\coborv : \rL(\X_{pq},A) \recht \rL(\X_{p+1,q},A) : (\coborv F)(X) &=
g_{10} \cdot F(\borv_0 X) + \sum_{j = 1}^{p+1} (-1)^j \; F(\borv_j X) \;
,
\end{align*}
when $X \in \X_{p,q+1}$ or $\X_{p+1,q}$ are as in Equation \eqref{eq.generalelement}. It is clear that we have in fact face operators
\begin{equation} \label{eq.faces}
\begin{split}
& \coborh_i : \rL(\X_{pq},A) \recht \rL(\X_{p,q+1}) \quad\text{for}\quad 0 \leq i \leq q+1 \quad\text{such that} \quad \coborh = \sum_{i=0}^{q+1} (-1)^i \coborh_i \; , \\
& \coborv_i : \rL(\X_{pq},A) \recht \rL(\X_{p+1,q},A)  \quad\text{for}\quad 0 \leq i \leq p+1 \quad\text{such that} \quad \coborv = \sum_{i=0}^{p+1} (-1)^i \coborv_i \; .
\end{split}
\end{equation}

For reasons that will become clear later, the elements of
$\rL(\X_{11},A)$ should be considered as $1$-cochains rather than
$0$-cochains. So, the total complex is defined as
\begin{equation} \label{eq.Kactotalcomplex}
C^n(A) = \bigoplus_{p+q = n+1} \rL(\X_{pq},A)
\end{equation}
and
\begin{equation} \label{eq.totalcomplexcoboundary}
\cobor : C^n(A) \recht C^{n+1}(A) : \cobor = \coborh  + \; \varepsilon
\coborv \quad\text{where}\;\; \varepsilon(F) = (-1)^q \quad\text{when}\quad
F \in \rL(\X_{pq},A) \; .
\end{equation}
We shall see in Proposition \ref{prop.nogwatzever} how
to define in a natural way $C^0 = A$ and $\cobor : C^0 \recht C^1$.

The following proposition is almost obvious.
\begin{proposition} \label{prop.groupofextensions}
The group of extensions of the matched pair $G_1,G_2 \subset G$ is precisely the second cohomology group of the total complex defined
by Equation \eqref{eq.Kactotalcomplex} with coefficients in the trivial $G$-module $\T$.
\end{proposition}
\begin{proof}
This is just a matter of making the right identifications, taking into account Remark \ref{rem.measurespace}. If $\cU : G_2 \times G_1 \times G_1 \recht \T$ and $\cV : G_2 \times G_2 \times G_1 \recht \T$ are measurable maps,
we define
$$\cU \Biggl(
\parbox[c]{6.5ex}{\begin{picture}(6.5,12)(1.5,0)
\put(2,2){\line(1,0){4}}
\put(2,6){\line(1,0){4}}
\put(2,10){\line(1,0){4}}
\put(2,2){\line(0,1){8}}
\put(6,2){\line(0,1){8}}
\put(6.4,3.4){\mbox{$h$}}
\put(6.4,7.4){\mbox{$g$}}
\put(2,10.4){\parbox[b]{4ex}{\centering $s$}}
\end{picture}}
\Biggr) = \cU(s,g,h)
\quad\text{and}\quad
\cV \Bigl(
\parbox[c]{10.5ex}{\begin{picture}(10.5,8)(1.5,0)
\put(2,2){\line(1,0){8}}
\put(2,6){\line(1,0){8}}
\put(2,2){\line(0,1){4}}
\put(6,2){\line(0,1){4}}
\put(10,2){\line(0,1){4}}
\put(10.2,3.4){\mbox{$g$}}
\put(2,6.4){\parbox[b]{4ex}{\centering $s$}}
\put(6,6.4){\parbox[b]{4ex}{\centering $t$}}
\end{picture}} \Bigr) = \cV(s,t,g) \; .
$$
Then, $(\cU,\cV) \in C^2(\T)$. The equation $\cobor (\cU,\cV) = 0$ precisely agrees with the three equations in \eqref{eq.cocycles}. Further, if $\cR : G_2 \times G_1 \recht \T$,
we define $\cR\bigl( \hspace{-2ex}\raisebox{0ex}[3.5ex][0ex]{\eenmatrix{s}{g}{{}}{{}}} \bigr) = \cR(s,g)$. The equation $\cobor \cR = (\cU,\cV)$ is equivalent with Equation
\eqref{eq.trivialcocycles}.
\end{proof}

\begin{remark} \label{rem.doublegroupoid}
The locally compact space $\X_{11}$ carries the structure of a {\it
  double groupoid}, \cite{brown1,brown2}. The {\it horizontal groupoid} $\Xh$ has unit space
  $\Xh^{(0)} = G_1$ which is embedded by $g \mapsto
  \raisebox{0ex}[3.5ex][0ex]{\eenmatrix{e}{g}{g}{e}}$. The source and range maps are defined by
$$\raisebox{0pt}[0.7cm][0.4cm]{}\source \Bigl(\eenmatrix{s}{g}{h}{t} \Bigr) = g \quad\text{and}\quad
\range \Bigl(\eenmatrix{s}{g}{h}{t} \Bigr) = h \; .$$
The composition is defined by
$$\eenmatrix{s}{g}{h}{t} \; \cdot\; \eenmatrix{s'}{h'}{g}{t'} \; = \;
\eenmatrix{ss'}{h'}{h}{tt'}$$
Analogously, the same space $\X_{11}$ carries a second groupoid
structure, the {\it vertical groupoid} $\Xv$, with unit space
$\Xv^{(0)}=G_2$. Composition is now defined by vertical composition of
squares. As such, $\X_{11}$ becomes a double groupoid: if $x,y,z,v$ are
such that $(x,y),(z,v) \in \Xh^{(2)}$ and $(x,z),(y,v) \in \Xv^{(2)}$,
then $(x \cdot\hor y, z \cdot\hor v) \in \Xv^{(2)}$, $(x \cdot\ver z,
y \cdot\ver v) \in \Xh^{(2)}$ and
$$(x \cdot\hor y) \cdot\ver (z \cdot\hor v) = (x \cdot\ver z)
\cdot\hor (y \cdot\ver v) \; .$$
This is obvious if you just look at the square
$$\parbox[c]{8ex}{\setlength{\unitlength}{1ex}
\begin{picture}(8,8)(0,0)
\put(0,0){\line(1,0){8}}
\put(0,4){\line(1,0){8}}
\put(0,8){\line(1,0){8}}
\put(0,0){\line(0,1){8}}
\put(4,0){\line(0,1){8}}
\put(8,0){\line(0,1){8}}
\put(1.3,5.3){\mbox{$x$}}
\put(5.3,1.3){\mbox{$v$}}
\put(5.3,5.3){\mbox{$y$}}
\put(1.3,1.3){\mbox{$z$}}
\end{picture}}$$

The bicomplex of Equation \eqref{eq.Kacbicomplex} can now be written
down analogously for double groupoids and gives a natural candidate for a
{\it double groupoid cohomology}.
\end{remark}

\section{The Kac exact sequence}

We still have a fixed matched pair $G_1,G_2 \subset G$ of locally
compact groups. We consider a Polish $G$-module $A$.

Looking at the bicomplex in Equation \eqref{eq.Kacbicomplex}, it is natural to add a row and a column and to write
\begin{equation}\label{eq.bigbicomplex}
\begin{CD}
\vdots @. \vdots @. \\
@A{\coborv}AA @A{\coborv}AA @. \\
\rL(\X_{10},A) @>{\coborh}>> \rL(\X_{11},A) @>{\coborh}>> \cdots\\
@A{\coborv}AA @A{\coborv}AA @. \\
\rL(\X_{00},A) @>{\coborh}>> \rL(\X_{01},A) @>{\coborh}>> \cdots
\end{CD}
\end{equation}

\begin{proposition} \label{prop.totalisG}
The cohomology of the total bicomplex of Equation \eqref{eq.bigbicomplex} is isomorphic with the measurable cohomology of $G$ with coefficients in the Polish $G$-module $A$. Moreover,
Equation \eqref{eq.isomorphism} defines an explicit cohomology
isomorphism. The inverse isomorphism is given explicitly in Remark \ref{rem.inverse}.
\end{proposition}
\begin{proof}
In order to prove that the cohomology of the total bicomplex of Equation \eqref{eq.bigbicomplex} is precisely the measurable cohomology of $G$, we consider the $G$-equivariant
bicomplex $(\rL(\Y_{pq},A))_{p,q \geq 0}$ and embed $A \hookrightarrow \rL(\Y_{00},A)$ as constant functions. We first prove that the total bicomplex completed with the
embedding $A \hookrightarrow \rL(\Y_{00},A)$ gives a free resolution of $A$.

By definition, every of the Polish $G$-modules $\rL(\Y_{pq},A) \cong \rL(G,\rL(\X_{pq},A))$ is free.

Consider now an arbitrary row in the bicomplex $(\rL(\Y_{pq},A))_{p,q
  \geq 0}$. Using the isomorphism of $\rL(\Y_{p,q+1},A)$ with $\rL(\Y_{pq} \times G_2,A)$ through
the identification of
$$ 
\parbox[c]{11ex}{\setlength{\unitlength}{1ex}\begin{picture}(9,6.5)(0,0)
\put(0,0){\line(1,0){9}}
\put(0,6){\line(1,0){9}}
\put(0,0){\line(0,1){6}}
\put(6,0){\line(0,1){6}}
\put(9,0){\line(0,1){6}}
\put(2,2.5){$X$}
\put(8.6,6.2){$x$}
\end{picture}} \in \Y_{p,q+1} \qquad\text{with}\qquad (X,p_2(x)) \in \Y_{pq} \times G_2$$
almost everywhere, we can write, for $F \in \rL(\Y_{pq},A)$, that $(\coborh F)(X,s) = \coborh(F(\cdot,s))(X) + (-1)^{q+1} F(X)$ almost everywhere. Hence, if $F \in \rL(\Y_{pq},A)$
and $\coborh F = 0$, we can use the Fubini theorem to take an $s \in G_2$ such that $0 = \coborh(F(\cdot,s))(X) + (-1)^{q+1} F(X)$ for almost all $X \in \Y_{pq}$. So,
the horizontal cohomology of the bicomplex $(\rL(\Y_{pq},A))_{p,q
  \geq 0}$ vanishes. It follows that the the total cohomology only
lives in the first column.

More precisely, this means that the total cohomology is the
cohomology of the complex $T^p = \{ F \in \rL(\Y_{p0},A) \mid \coborh F = 0 \}$ with $\coborv$ as a coboundary operator.
Completing with $A \hookrightarrow \rL(\Y_{00},A)$,
we claim that we precisely get the standard resolution for the measurable cohomology of $G_1$ with coefficients in $A$. Observe that $\rL(\Y_{p0},A) \cong \rL(G_2 \times G_1^{p+1},A)$ through
the identification of the column vector $(x_0,\ldots,x_p)$ in $\Y_{p0}$ with $(p_2(x_0),p_1(x_0),\ldots,p_1(x_p)) \in G_2 \times G_1^{p+1}$ almost everywhere. It is then
easy to check that $T^p \cong \rL(G_1^{p+1},A)$ and that this isomorphism intertwines $\coborv$ with the usual group coboundary operator. This proves our claim. So, we have
proven that the total bicomplex of $(\rL(\Y_{pq},A))_{p,q \geq 0}$, completed with $A \hookrightarrow \rL(\Y_{00},A)$ gives a free resolution of $A$. Hence, the measurable
group cohomology $H(G,A)$ is given as the total cohomology of the bicomplex $(\rL(\Y_{pq},A)^G)_{p,q \geq 0}$, which is precisely $(\rL(\X_{pq},A))_{p,q \geq 0}$.

Denote the total bicomplex of \eqref{eq.bigbicomplex} by $D(A) = (D^n(A))_{n \geq 0}$. We have proved that $H(D,A) \cong H(G,A)$. But there is more. Since the cohomology
$H(D,A)$ turns short exact sequences of coefficient modules in a natural way into short exact cohomology sequences, the cohomology theory $H(D,A)$ satisfies the Buchsbaum criterion.
We can calculate $H(G,A)$ using the
complex $(\funborel(G^n,A))_{n \geq 0}$, where $\funborel$ denotes all
Borel measurable functions and where the coboundary maps are defined
in the Preliminaries. Whenever now $I:(\funborel(G^n,A)) \recht (D^n(A))$ is a
natural cochain transformation, which is the identity on $\funborel(G^0,A) = A =
D^0(A)$, we can conclude that this cochain transformation is a cohomology
isomorphism. Such a cochain transformation can be written as
\begin{equation} \label{eq.isomorphism}
I : \funborel(G^n,A) \recht D^n(A) : I(F)(X) = \sum_{\text{path in}\;
  X} \text{Sign} (\text{path}) \; F(\text{path}) \; ,
\end{equation}
where $X \in \X_{pq}$, $p+q=n$ and a path in e.g. $X \in \X_{23}$ is a thing like

\begin{center}
\raisebox{0cm}[1cm][1cm]{\hspace{20pt}\parbox[c]{120pt}{\setlength{\unitlength}{1mm}
\begin{picture}(0,0)(0,0)
\put(3,21.5){\mbox{$s$}}
\put(14,11.5){\mbox{$t$}}
\put(24,11.5){\mbox{$r$}}
\put(12.5,16){\mbox{$g$}}
\put(32.5,6){\mbox{$h$}}
\end{picture}
\includegraphics[width=3.2cm]{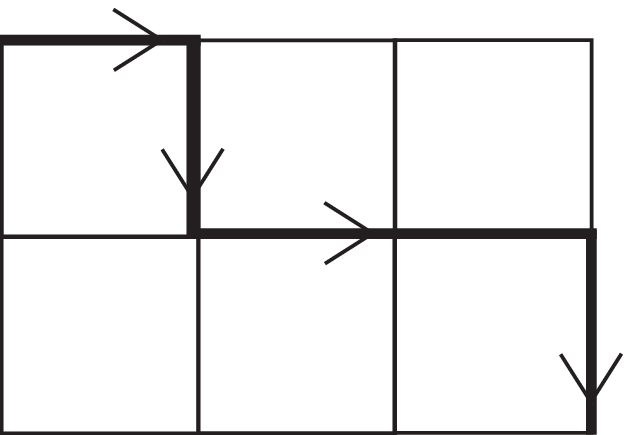}}
\parbox[c]{270pt}{For this path $F(\text{path}) = F(s,g,t,r,h)$ and
$\text{Sign}(\text{path})$ equals \\ $(-1)^{\text{number of squares above the
    path}}$, i.e., $(-1)^2=1$.}}
\end{center}

More formally, a path in $X \in \X_{pq}$ is a path that starts in the
top left corner of $X$, goes either down or right and ends, after
$p+q=n$ steps in the bottom right corner of $X$. We evaluate $F
\in \funborel(G^n,A)$ on the edges along the
path of length $n$. The sign of the path
is defined as $-1$ to the power the number of squares that are
above the path. An elementary computation learns that $I$ is a cochain transformation.
\end{proof}

\begin{remark} \label{rem.inverse}
We can also write a natural cochain transformation
$I' : D^n(A) \recht \rL(G^n,A)$, yielding the inverse isomorphism
$H(D(A)) \recht H(G,A)$. Define, the measure class isomorphism $G^n
\recht \X_{nn} : (x_1,\ldots,x_n) \mapsto X(x_1,\ldots,x_n)$ where
\begin{equation*}
X(x_1) = \kadermetx{\, x_1} \quad , \quad
X(x_1,x_2) =
\parbox[c]{8ex}{\setlength{\unitlength}{1ex}
\begin{picture}(8,8)(0,0)
\put(0,0){\line(1,0){8}}
\put(0,4){\line(1,0){8}}
\put(0,8){\line(1,0){8}}
\put(0,0){\line(0,1){8}}
\put(4,0){\line(0,1){8}}
\put(8,0){\line(0,1){8}}
\put(1,5.5){\mbox{$x_1$}}
\put(5,1.5){\mbox{$x_2$}}
\end{picture}}
\end{equation*}
and where $X(x_1,\ldots,x_n) \in \X_{nn}$ is defined analogously by
putting $x_1,\ldots,x_n$ in boxes along the diagonal. Further, we
define $P_i : \X_{nn} \recht \X_{n-i,i}$ where $P_i X$ is the lower
left corner of $X$, i.e.\ consisting of the $n-i$ final lines and $i$
first columns of $X$. We can now define
$$I' : D^n(A) \recht \rL(G^n,A) : I'(F)(x_1,\ldots,x_n) = \sum_{i=0}^n
F\bigl(P_i(X(x_1,\ldots,x_n))\bigr) \; .$$
One verifies that $I'$ is a cochain transformation.
\end{remark}

Observe that $\X_{0n} = G_2^n$ and that $\coborh : \rL(\X_{0n},A) \recht
\rL(\X_{0,n+1},A)$ coincides with the usual coboundary operator $\cobor
: \rL(G_2^n, A) \recht \rL(G_2^{n+1},A)$. Making the same observation for the
first column of the bicomplex of Equation \eqref{eq.bigbicomplex}, we
obtain a natural cochain transformation
\begin{equation} \label{eq.chaintransformation}
J : (D^n(A))_n \recht (K^n(A))_n \quad\text{where}\quad
D^n(A)=\bigoplus_{p+q=n} \rL(\X_{pq},A) \quad\text{and}\quad
K^n(A)=\rL(G_1^n,A) \oplus \rL(G_2^n,A)
\end{equation}
and where we just have to explain that, for $n=0$, we take $A \ni a
\mapsto a \oplus a \in A \oplus A$.

The following definition may seem a bit pedantic, but implies
in a natural way how to define the Kac $0$-cohomology (see Proposition
\ref{prop.nogwatzever} for the link with the cohomology of the cochain
complex $(C^n(A))$ defined in Equations \eqref{eq.Kactotalcomplex} and
\eqref{eq.totalcomplexcoboundary}).

\begin{definition}\label{def.Kacchaincomplex}
The Kac cohomology $H(\matched,A)$ of the matched pair ($\matched$) $G_1,G_2 \subset G$ with
coefficients in the Polish $G$-module $A$ is defined to be the
cohomology of the mapping
cone of the natural cochain transformation $J$ given by Equation
\eqref{eq.chaintransformation}.
\end{definition}

We recall that, by definition, the mapping cone $(M^n)_{n \geq -1}$ of
the cochain transformation $J$,
is defined by the formula
$$M^n(A) = D^{n+1}(A) \oplus K^n(A) \quad\text{and}\quad \cobor(F,G) = (\cobor
F, -\cobor G + J F) \; .$$

We now explain how to complete the total bicomplex of \eqref{eq.Kacbicomplex} in order to obtain the Kac cohomology.

\begin{proposition} \label{prop.nogwatzever}
Define, for $n \geq 1$, $C^n(A)$ as in Equation
\eqref{eq.Kactotalcomplex}. Define $C^0(A) = A$, $C^{-1}(A) = 0$. Define
$$\cobor : C^0(A) \recht C^1(A) : (\cobor
a)\Bigl(\eenmatrixmetx{s}{g}{h}{t}{x} \Bigr) = x \cdot a - s \cdot a -
h \cdot a + a \; .$$
This coboundary compiles with Equation
\eqref{eq.totalcomplexcoboundary} to a cochain complex $(C^n(A))_{n \geq -1}$.
The inclusions $C^n(A) \hookrightarrow D^{n+1}(A)$, together with the
map
$$C^0(A) \recht M^0(A) = \rL(\X_{10},A) \oplus \rL(\X_{01},A) \oplus A
\oplus A :
a \mapsto \coborv a \oplus \coborh a \oplus a \oplus a$$
define a natural cochain transformation $(C^n(A)) \recht (M^n(A))$
which is a cohomology isomorphism.

In particular, we conclude from Proposition
\ref{prop.groupofextensions} that the group of extensions of the
matched pair $G_1,G_2 \subset G$ is precisely the Kac $2$-cohomology
$H^2(\matched,\T)$ with coefficients in the trivial $G$-module $\T$.
\end{proposition}
\begin{proof}
For $n \geq 2$, $J : D^n(A) \recht K^n(A)$ is surjective, with kernel
the image of $C^{n-1}(A)$. This shows that the cochain transformation
$(C^n(A)) \recht (M^n(A))$ is an isomorphism in $n$-cohomology for $n
\geq 2$. For $n=0,1$, the same follows by an explicit verification.
\end{proof}

The mapping cone of a cochain transformation is made for getting long
exact sequences. So, the following is an immediate consequence of
Proposition \ref{prop.totalisG}.

\begin{corollary} \label{cor.Kacexact}
The Kac cohomology $H(\matched,A)$ of the matched pair $G_1,G_2 \subset G$ with
coefficients in the Polish $G$-module $A$ satisfies the long exact
sequence
\begin{align*}
0 & \recht A^G \recht A^{G_1} \oplus A^{G_2} \recht H^0(\matched,A) \recht H^1(G,A)
\recht H^1(G_1,A) \oplus H^1(G_2,A) \recht H^1(\matched,A) \\ & \recht H^2(G,A)
\recht H^2(G_1,A) \oplus H^2(G_2,A) \recht H^2(\matched,A) \recht H^3(G,A)
\recht H^3(G_1,A) \oplus H^3(G_2,A) \recht \cdots
\end{align*}
Recall here that $H^2(\matched,\T)$ is the group of extensions of the
matched pair, by Propositions \ref{prop.groupofextensions} and \ref{prop.nogwatzever}.
\end{corollary}

\begin{remark}
Once we have, for a good class of locally compact groups $G$ and
Polish $G$-modules $A$, a natural way to write a cochain complex
$(E^n(G,A))$ whose cohomology is the cohomology of $G$, we can expect
that the Kac cohomology of a matched pair $G_1,G_2 \subset G$ is the
cohomology of the mapping cone of the cochain transformation
$(E^n(G,A))_n\recht (E^n(G_1,A) \oplus E^n(G_2,A))_n$. Concrete
applications of this principle can be found below (see Proposition \ref{prop.mappingcone}).
\end{remark}

If $A$ is a Polish $G$-module, we define $\rL(G,A)$ as a Polish
$G$-module by writing $(x \cdot F)(y) = x \cdot F(x^{-1}y)$ for $F \in
\rL(G,A)$ and $x,y \in G$. Sending $a \in A$ to the constant function
$a$, we get an embedding $A \hookrightarrow \rL(G,A)$ of $A$ as a
closed submodule of $\rL(G,A)$. As a consequence of Corollary
\ref{cor.Kacexact}, we get the following result.

\begin{corollary} \label{cor.kacefface}
If $G_1,G_2 \subset G$ is a matched pair, then $H^n(\matched,\rL(G,A)) = 0$
for all $n \geq 1$ and all Polish $G$-modules $A$. In particular, the
Kac cohomology satisfies the Buchsbaum criterion.
\end{corollary}
\begin{proof}
Since the group cohomology with coefficients in $\rL(G,A)$ vanishes,
this is an immediate consequence of the Kac exact sequence.
\end{proof}

It is now possible to interpret the Kac cohomology of the matched pair $G_1,G_2 \subset G$ with coefficients in the Polish $G$-module $A$ as an ordinary
group cohomology with coefficients in a well chosen module.

\begin{corollary} \label{cor.Kacasgroup}
Let $G_1,G_2 \subset G$ be a matched pair and $A$ a Polish $G$-module. Then the Kac cohomology $H(\matched,A)$ is the measurable group cohomology of $G$
with coefficients in the Polish $G$-module
$$\frac{\rL(G/G_1 \sqcup G/G_2,A)}{A} \; .$$
\end{corollary}
\begin{proof}
By Proposition \ref{prop.nogwatzever},
the Kac cohomology is defined by the total bicomplex of $(\rL(\X_{pq},A))_{p,q \geq 1}$ completed with $A \cong \rL(\X_{00},A) \recht \rL(\X_{11},A) : a \mapsto \cobor(a) =
(\coborv_1 - \coborv_0)(\coborh_1 - \coborh_0) a$. Consider now the $G$-equivariant analogue, which is the total bicomplex of $(\rL(\Y_{pq},A))_{p,q \geq 1}$ completed with
$\rL(\Y_{00},A) \recht \rL(\Y_{11},A) : F \mapsto \cobor F = (\coborv_1 - \coborv_0)(\coborh_1 - \coborh_0) F$. Denote the complex obtained as such by $(R^n)$.

We claim that
$$\frac{\rL(G/G_1 \sqcup G/G_2,A)}{A} \overset{\te}{\recht} R^0 \recht R^1 \recht \cdots \qquad\text{where}\quad \te(H_1 \oplus H_2) = H_1 - H_2$$
is a free resolution. The proof of this claim will complete the proof of the corollary.

By definition, every of the Polish $G$-modules $R^n$ is free. Moreover, the exactness of $R^{n-1} \recht R^{n} \recht R^{n+1}$ for $n \geq 1$, follows
from Corollary \ref{cor.kacefface} stating that
$H^n(\matched,\rL(G,A)) = 0$. Also, $\te$ is injective. If $H_1 \in
\rL(G/G_1,A)$ and $H_2 \in \rL(G/G_2,A)$ are such that $\te(H_1 \oplus
H_2) = 0$, then $H_1 = H_2 = F$, where $F \in \rL(G,A)$ is invariant
under translations by both $G_1$ and $G_2$. Hence, $F$ is a constant
function, proving the injectivity of $\te$.

We have to show that the kernel of $\cobor :
\rL(\Y_{00},A) \recht \rL(\Y_{11},A)$ is the image of $\te$. It is
immediate the image of $\te$ is included in the kernel of $\cobor$.
Let $F \in \rL(\Y_{00},A)$ and $\cobor F = 0$. This means that
$$F(xy) - F(xq_2(y)) - F(xp_1(y)) + F(x) = 0 \qquad\text{for almost all}\quad (x,y) \in G \times G \; .$$
Using the Fubini theorem, we fix $x \in G$ such that the previous equation holds for almost all $y \in G$. Define $H_1 \in \rL(G/G_1,A)$ and $H_2 \in \rL(G/G_2,A)$ by
the formulas $H_1(y) = F(xq_2(y))$ and $H_2(y) = F(x) - F(xp_1(y))$. Then, $F = \te(H_1 \oplus H_2)$.
\end{proof}

\begin{remark}
Observe that it is by now clear that for any Polish $G$-module $A$, we
have $H(G,\rL(G/G_i,A)) = H(G_i,A)$. Hence, the Kac exact sequence is
exactly the long exact cohomology sequence that corresponds to the
exact sequence of Polish $G$-modules
$$0 \recht A \recht \rL(G/G_1 \sqcup G/G_2,A) \recht \frac{\rL(G/G_1
  \sqcup G/G_2,A)}{A} \recht 0 \; .$$ 
\end{remark}

\section{Pentagonal cohomology and isomorphism with Kac cohomology}

In Definition \ref{def.pentagonalcohomology}, we defined the $2$-cohomology group associated with a pentagonal transformation.
Recall Equation \eqref{eq.faces}, where we defined face operators $\coborh_i : \rL(\X_{pq},A) \recht \rL(\X_{p,q+1},A)$ and
$\coborv_j : \rL(\X_{pq},A) \recht \rL(\X_{p+1,q},A)$ for $0 \leq i \leq q+1$ and $0 \leq j \leq p+1$. For convenience of notation, we write moreover
$\coborh_{q+2} = \coborh_{q+1}$ and $\coborv_{-1} = \coborv_0$ on $\rL(\X_{pq},A)$. Define then, for $n \geq 1$,
$$\coborpent : \rL(\X_{nn},A) \recht \rL(\X_{n+1,n+1},A) : \coborpent = \sum_{i=0}^{n+2} (-1)^i \coborh_i \coborv_{i-1} \; .$$
In Remark \ref{rem.inverse}, we defined the measure class isomorphism $G^n \recht \X_{nn}$. Using this, we can identify the usual complex for the measurable
cohomology of $G$ with the complex defined by
$$\cobor\group : \rL(\X_{nn},A) \recht \rL(\X_{n+1,n+1},A) : \cobor\group = \sum_{i=0}^{n+1} (-1)^i \coborh_i \coborv_i \; .$$
It is then clear that the injective map $\coborh_{n+1}\coborv_0 : \rL(\X_{nn},A) \recht \rL(\X_{n+1,n+1},A)$ intertwines $\coborpent$ and $\cobor\group$. Hence, $\coborpent$
turns $E^n(A):=\rL(\X_{nn},A)$ for $n \geq 1$ into a cochain complex. Remark that the $2$-cohomology of this complex is, by definition, the $2$-cohomology associated
with the pentagonal transformation.

We can complete the cochain complex $(E^n(A))_{n \geq 1}$, defining $E^0(A) = A \oplus A$ with
$$\coborpent : A \oplus A \recht \rL(\X_{11},A) : \coborpent(a,b) = \cobor(a) + \coborh_0 \coborv_1 (b) \; ,$$
where $\cobor : A \recht \rL(\X_{11},A)$ is as in Proposition
\ref{prop.nogwatzever}.

\begin{definition}
Given a matched pair $G_1,G_2 \subset G$, we define the
\emph{pentagonal cohomology} $H(\pent,A)$ with coefficients in the
Polish $G$-module $A$ as the cohomology of the cochain complex
$(E^n(A))_{n \geq 0}$ defined above.
\end{definition}

Recall the Kac cochain complex $(C^n(A))$ that we defined in Equation
\eqref{eq.Kactotalcomplex} and Proposition
\ref{prop.nogwatzever}. Using Remark \ref{rem.inverse} and the cochain
transformation $I' : D^n(A) \recht \rL(G^n,A) \cong \rL(\X_{nn},A)$
defined there, we get a cochain transformation $T : C^n(A) \recht E^n(A)$, which
is defined such that for $n \geq 1$, the diagram
\begin{equation*}
\begin{CD}
C^n(A) @>>> D^{n+1}(A) \\
@V{T}VV @VV{I'}V \\
E^n(A) @>{\coborh_{n+1}\coborv_0}>> \rL(\X_{n+1,n+1},A) \cong  \rL(G^{n+1},A)
\end{CD}
\end{equation*}
commutes. For $n=0$, we put $T : C^0(A) = A \recht A \oplus A = E^0(A) : a \mapsto a \oplus 0$.

We prove now that $T$ is a cohomology
isomorphism. Hence, the pentagonal and Kac cohomologies are
isomorphic.

\begin{proposition} \label{prop.iso}
Let $G_1,G_2 \subset G$ be a matched pair. The
cochain transformation $T$ is a cohomology isomorphism between Kac
cohomology and pentagonal cohomology.
\end{proposition}
\begin{proof}
As in Equation \eqref{eq.equivariant}, we can write equivariant face
operators
\begin{align*}
& \coborh_i : \rL(\Y_{pq},A) \recht \rL(\Y_{p,q+1},A)
\quad\text{for}\quad 0 \leq i \leq q+1 \qquad\text{and} \\
&\coborv_j : \rL(\Y_{pq},A) \recht \rL(\Y_{p+1,q},A)
\quad\text{for}\quad 0 \leq j \leq p+1 \; .
\end{align*}
For convenience of notation, we put moreover $\coborh_{q+2} =
\coborh_{q+1}$ and $\coborv_{-1} = \coborv_0$ on
$\rL(\Y_{pq},A)$. Using this notation, we get the $G$-equivariant
cochain complex $(P^n(A))$ such that for $n \geq 1$, $P^n(A) =
\rL(\Y_{nn},A)$ and
$$\coborpent : \rL(\Y_{nn},A) \recht \rL(\Y_{n+1,n+1},A) : \coborpent
= \sum_{i=0}^{n+2} (-1)^i \coborh_i \coborv_{i-1} \; .$$
Also, write $P^0 = \rL(\Y_{00},A) \oplus \rL(\Y_{00},A)$ and
$$\coborpent : \rL(\Y_{00},A) \oplus \rL(\Y_{00},A) \recht
\rL(\Y_{11},A) : \coborpent(\al \oplus \beta) = (\coborh_0 -
\coborh_1)(\coborv_0 - \coborv_1) \al + \coborh_0 \coborv_1 \be \; .$$
Then, the cochain complex $(E^n(A))$ of the pentagonal cohomology
consists precisely of the fixed points under $G$ in $(P^n(A))$.

We now complete the above $G$-equivariant complex by $\frac{\rL(G/G_1 \,\sqcup\, G/G_2,A)}{A}$ in order to get a resolution.
Using Corollary \ref{cor.Kacasgroup}, it is sufficient
to prove that
$$\frac{\rL(G/G_1 \sqcup G/G_2,A)}{A} \overset{\pi}{\recht} P^0(A) \recht P^1(A) \recht \cdots
\qquad\text{with}\quad \pi(H_1 \oplus H_2) = (H_1 - H_2) \oplus 0$$
is a free resolution of $\frac{\rL(G/G_1 \,\sqcup\, G/G_2,A)}{A}$. Indeed, we then know that the pentagonal cohomology is isomorphic with the Kac cohomology and that
the natural cochain transformation $T$ induces a cohomology isomorphism.

From
the lemma following this proposition, we know that the sequence $P^0(A) \recht P^1(A) \recht \cdots$ is exact. It remains to prove that if $F \oplus G \in P^0(A)$ and
$\cobor F + \coborh_0 \coborv_1 G = 0$, then $F \oplus G$ belongs to the image of $\pi$. From the proof of Corollary \ref{cor.Kacasgroup}, it follows that it is sufficient to
prove that $G=0$. But, $(\coborv_1 - \coborv_0)(\coborh_1 - \coborh_0)F + \coborh_0 \coborv_1 G = 0$. Apply $\coborh = \coborh_0 - \coborh_1 + \coborh_2$ to both sides of this
last equation. We find that $0 =  \coborh_2 \coborh_0 \coborv_1 G$ and hence, $G=0$.
\end{proof}

\begin{lemma}
Let $G_1,G_2 \subset G$ be a matched pair and $A$ a Polish $G$-module. The sequence $$P^0(A) \recht P^1(A) \recht \cdots$$ constructed in the proof of the
previous proposition, is exact.
\end{lemma}
\begin{proof}
We temporarily consider some cohomology
which is intermediate between the pentagonal cohomology and the group
cohomology of $G$. For $n \geq 0$, we define
$$\coborhelp : \rL(\Y_{n+1,n},A) \recht \rL(\Y_{n+2,n+1},A) :
\coborhelp = \sum_{i=0}^{n+2} (-1)^i \coborh_i \coborv_i \; .$$
The injective maps
$$\rL(\Y_{nn},A) \overset{\coborv_0}{\longrightarrow} \rL(\Y_{n+1,n},A)
\overset{\coborh_{n+1}}{\longrightarrow} \rL(\Y_{n+1,n+1},A)$$
intertwine $\coborpent,\coborhelp$ and $\coborgroup$ for $n \geq 1$.

Suppose $n \geq 1$ and write $\coborhelp^n$ for a while. We identify
$\rL(\Y_{n+2,n+1},A)$ with $\rL(G \times \Y_{n+1,n},A)$ through the
identification of
$$(x,X) \in G \times \Y_{n+1,n} \qquad\text{with}\qquad
\parbox[c]{11ex}{\setlength{\unitlength}{1ex}\begin{picture}(11,8)(0,0)
\put(2,0){\line(1,0){8}}
\put(2,6){\line(1,0){8}}
\put(2,8){\line(1,0){8}}
\put(2,0){\line(0,1){8}}
\put(4,0){\line(0,1){8}}
\put(10,0){\line(0,1){8}}
\put(0.7,8){$x$}
\put(6,2.5){$X$}
\end{picture}} \in \Y_{n+2,n+1}$$
almost everywhere. Observe that
$$(\coborhelp^n F)(x,X) = F(X) -
\coborhelp^{n-1}\bigl(F(x,\cdot)\bigr)(X) \quad\text{for almost all}\;
(x,X) \in G \times \Y_{n+1,n} \; .$$
If now $F \in \rL(\Y_{n+1,n},A)$ and $\coborhelp^n F = 0$, we take a
Borel measurable representative for $F$ that we still denote by
$F$. Then,
$$F(X) = \coborhelp^{n-1}\bigl(F(x,\cdot)\bigr)(X) \quad\text{for almost all}\;
(x,X) \in G \times \Y_{n+1,n} \; .$$
By the Fubini theorem, we can take an $x \in G$ such that the previous
equation holds for almost all $X \in \Y_{n+1,n}$. So, we conclude that
$F$ is a coboundary.

As such, we did not only prove that the
cohomology of $\coborhelp$ is trivial, but we also proved the
following: if $\te : A \recht B$ is a continuous $G$-equivariant homomorphism between
the Polish $G$-modules $A$ and $B$ and if $F \in \rL(\Y_{n+1,n},A)$
such that $\coborhelp(\te(F)) = 0$, then there exists $K \in
\rL(\Y_{n,n-1},A)$ such that $\te(F) = \te(\coborhelp K)$.

Consider
\begin{equation}\label{eq.negative}
\coborhelp : \rL(G_1,A) \recht \rL(\Y_{10},A) : (\coborhelp
F)\bigl(\parbox[c]{3ex}{\setlength{\unitlength}{1ex}\begin{picture}(2,3)
\put(0.2,0){\line(0,1){3}}
\put(0.9,-0.5){$y$}
\put(0.9,3){$x$}
\end{picture}}\bigr) = F(p_1(x)) - F(p_1(y)) \; .
\end{equation}
If now $F \in
\rL(\Y_{10},A)$ and $\coborhelp F = 0$, there exists $K \in
\rL(G_1,A)$ such that $F = \coborhelp K$. The same kind of statement
involving $\te : A \recht B$ as above also holds.

We now prove the exactness of the sequence $P^0(A) \recht P^1(A) \recht \cdots$. Let $n \geq 1$ and
$\al \in \rL(\Y_{nn},A)$ such that $\coborpent \al = 0$. Then,
$\coborv_0 \al \in \rL(\Y_{n+1,n},A)$ and $\coborhelp \coborv_0 \al =
0$. Using the results above, we can take $\be \in \rL(\Y_{n,n-1},A)$
such that $\coborv_0 \al = \coborhelp \be$.

Suppose now first that $n \geq 2$. Applying $\coborv_0 -
\coborv_1$ to both sides of this equation yields
\begin{equation} \label{eq.torewrite}
0 = \sum_{i=1}^{n +1} (-1)^i \coborh_i \coborv_{i+1}
(\coborv_0-\coborv_1) \be \; .
\end{equation}
We now use the identification of $\rL(\Y_{n+1,n-1},A)$ with
$\rL(\Y_{n-1,n-2}, \rL(\Y_{10},A))$ by identifying
$$(X,Y) \in \Y_{n-1,n-2} \times \Y_{10} \qquad\text{with}\qquad
\parbox[c]{11ex}{\setlength{\unitlength}{1ex}\begin{picture}(11,8)(0,0)
\put(2,0){\line(1,0){8}}
\put(2,5){\line(1,0){8}}
\put(2,6.5){\line(1,0){8}}
\put(2,8){\line(1,0){8}}
\put(2,0){\line(0,1){8}}
\put(3.5,0){\line(0,1){8}}
\put(10,0){\line(0,1){8}}
{\linethickness{2pt}\put(2,6.3){\line(0,1){1.9}}
\put(3.5,0){\line(1,0){6.5}}
\put(3.5,0){\line(0,1){5}}
\put(3.5,5){\line(1,0){6.5}}
\put(10,0){\line(0,1){5}}
}
\put(6,1.5){$X$}
\put(-0.2,6.5){$Y$}
\end{picture}} \in \Y_{n+1,n-1}$$
almost everywhere. We also use the identification of
$\rL(\Y_{n,n-1},A)$ with $\rL(\Y_{n-1,n-2}, \rL(\Y_{00},A))$ as
above. We consider $\te := \coborv_0 - \coborv_1 : \rL(\Y_{00},A)
\recht \rL(\Y_{10},A)$ as a morphism of Polish modules. Using all
these identifications, Equation \eqref{eq.torewrite} becomes
$\coborhelp (\te(\be)) = 0$. From the results in the beginning of the
proof, we find $\ga \in \rL(\Y_{n-2,n-3},\rL(\Y_{00},A))$ such that
$\te(\be) = \te(\coborhelp \ga)$. Re-identifying, we have found $\ga
\in \rL(\Y_{n-1,n-2},A)$ such that
$$(\coborv_0 - \coborv_1) \be = \sum_{i=1}^n (-1)^i \coborh_i
\coborv_{i+1} (\coborv_0 - \coborv_1) \ga$$
and hence,
$$(\coborv_0 - \coborv_1) \be = (\coborv_0 - \coborv_1) \coborhelp \ga
\; .$$
Observe that the argument works well for $n=2$ by using Equation \eqref{eq.negative}.

Write $\betatil = \beta - \coborhelp \ga$. Then $\coborv_0 \al =
\coborhelp \be = \coborhelp \betatil$ and moreover $(\coborv_0 -
\coborv_1) \betatil = 0$. But then, there is a unique $\mu \in
\rL(\Y_{n-1,n-1},A)$ such that $\betatil = \coborv_0 \mu$. We conclude
that $\coborv_0 \al = \coborhelp \coborv_0 \mu = \coborv_0 \coborpent
\mu$ and hence, $\al = \coborpent \mu$.

Finally, consider the case $n=1$. So, we have $\al \in \rL(\Y_{11},A)$
such that $\coborpent \al = 0$. The beginning of the argument above
works well and we find $\be \in \rL(\Y_{10},A)$ such that
$$\coborv_0 \al = (\coborh_0 \coborv_0 - \coborh_1\coborv_1 +
\coborh_1 \coborv_2) \be \; .$$
Apply $\coborv_0 - \coborv_1$ to both sides and conclude that
$\coborh_1 (\coborv_0 - \coborv_1) (\coborv_1-\coborv_2)\be =
0$. Hence, $(\coborv_2 - \coborv_3)(\coborv_0 - \coborv_1) \be =
0$. So, we can take $\ga \in \rL(\Y_{10},A)$ such that
\begin{equation} \label{eq.hulpstuk}
(\coborv_0 -
\coborv_1) \be = \coborv_2 \ga \; .
\end{equation}
It follows that $\coborv \be = \coborv_2(\ga + \be)$ and hence
$\coborv \coborv_2(\ga + \be) = 0$. This can be rewritten as
$\coborv_3(\coborv_0 - \coborv_1)(\ga + \be) = 0$ and hence,
$(\coborv_0 - \coborv_1)(\ga + \be) = 0$. Combining this last equation
with Equation \eqref{eq.hulpstuk}, we find that $\coborv \ga = 0$. It
is by now easy to find $\mu \in \rL(\Y_{00},A)$ such that $\ga =
\coborv \mu$ and hence, $\ga = (\coborv_0 - \coborv_1)\mu$. Write
$\betatil = \beta - \coborv_1\mu$. Then,
$$(\coborv_0 - \coborv_1) \betatil = (\coborv_0 - \coborv_1)(\be +
\ga) - (\coborv_0 - \coborv_1)(\ga + \coborv_1 \mu) = 0 - (\coborv_0 -
\coborv_1)\coborv_0 \mu = 0 \; .$$
So, we can take $\eta \in \rL(\Y_{00},A)$ such that $\betatil =
\coborv_0 \eta$. Finally,
$$\coborhelp \betatil = \coborhelp \be - (\coborh_0 \coborv_0 -
\coborh_1(\coborv_1 - \coborv_2)) \coborv_1 \mu = \coborv_0(\al -
\coborh_0\coborv_1 \mu) \; .$$
Hence,
$$(\coborh_0 \coborv_0 - \coborh_1 \coborv_0 + \coborh_1
\coborv_1)(\eta) = \al - \coborh_0 \coborv_1 \mu$$
and so, $\al = d(\eta) + \coborh_0\coborv_1 (\eta + \mu) =
\coborpent(\eta \oplus (\eta+\mu))$.
\end{proof}

\section{Interpretation of $\mathbf{1}$-cohomology}

Fix a matched pair $G_1,G_2 \subset G$. We give a natural
interpretation of the Kac $1$-cohomology $H^1(\matched,\T)$ with
coefficients in the trivial $G$-module $\T$.

As we explained in Subsection \ref{subsec.bicros}, the elements of
$H^2(\matched,\T)$ can be interpreted as extensions $(M,\de)$
$$
e \longrightarrow (\rL^\infty(G_1),\de_1) \longrightarrow (M,\de)
\longrightarrow (\cL(G_2),\deh_2) \longrightarrow e
$$
in the category of locally compact quantum groups.
In particular, the element $0 \in H^2(\matched,\T)$ gives rise to an
extension $(M,\de)$ that we describe explicitly as follows. It is
called the bicrossed product of $G_1,G_2$.

The von Neumann algebra $M$ is the crossed product $M = G_2 \ltimes
\rL^\infty(G/G_2)$ that we realize explicitly as follows.
Identifying $\rL^\infty(G/G_2)$ with $\rL^\infty(G_1)$,
we consider the Hilbert space $H = L^2(G_2 \times G_1)$ and write
$$\pi_1 : \rL^\infty(G_1) \recht \B(H) : (\pi_1(F) \xi)(s,g) =
F(p_1(sg)) \xi(s,g) \quad\text{and}\quad \lambda : \cL(G_2) \recht \B(H) :
\lambda(a) = a \ot 1 \; .$$
The von Neumann algebra $M$ is generated by $\pi_1(\rL^\infty(G_1))$ and
$\lambda(\cL(G_2))$ in $\B(H)$.
The coproduct $\de : M \recht M \ot M$ on $M$ is then given by
\begin{align*}
&\de \pi_1 = (\pi_1 \ot \pi_1) \de_1 \quad\text{where}\quad
\de_1(F)(g,h) = F(gh) \quad\text{for all}\quad g,h \in G_1 \; , \\
&\de(\lambda_s \ot 1) = \bigl((\lambda_s \ot 1 \ot 1)(\pi_1 \ot
\io)(X_s)\bigr) \ot 1 \\ &\qquad\qquad\qquad\qquad
\text{where}\quad X_s \in \rL^\infty(G_1) \ot \cL(G_2)
\quad\text{is given by}\quad X_s(g) = \lambda_{p_2(sg)} \; .
\end{align*}

Since $M$ is a crossed product, there is a dual coaction
$$\te : M \recht \cL(G_2) \ot M : \begin{cases} \te(\pi_1(F)) = 1 \ot
  \pi_1(F) \; , \\ \te(\lambda_s \ot 1) = \lambda_s \ot (\lambda_s \ot
  1) \; .\end{cases}$$
If we write somehow formally the morphism $\pi_2 : M \recht \cL(G_2)$,
we can consider that $\te = (\pi_2 \ot \io)\de$. This is explained in
detail after Proposition 3.1 in \cite{VV}.

\begin{definition}
A (left) coideal $I$ of a locally compact quantum group $(M,\de)$ is a von Neumann
subalgebra $I \subset M$ satisfying $\de(I) \subset M \ot I$.
\end{definition}
We only work with left coideals and so, we leave out \lq left\rq\ from
now on.

We shall give a bijective correspondence between the Kac
$1$-cohomology $H^1(\matched,\T)$ and a special class of coideals in
the bicrossed product $(M,\de)$. If $I$ is a coideal, then $\te(I)
\subset \cL(G_2) \ot I$. So, $(\cL(G_2),\deh_2)$ coacts on $I$. As
such, $I$ is a $\cL(G_2)$-comodule.

\begin{definition}
We say that a coideal $I$ of the bicrossed product $(M,\de)$ is
\emph{full}, if for every $s \in G_2$,
$$E_s := \{ x \in I \mid \te(x) = \lambda_s \ot x \}$$
is a one-dimensional subspace of $I$.
\end{definition}

We then have the following result.

\begin{proposition} \label{prop.interprete}
Let $I$ be a coideal of $(M,\de)$. Then, the following are equivalent.
\begin{enumerate}
\item $I$ is full.
\item There exists an isomorphism $\mu : \cL(G_2) \recht I$ of von
  Neumann algebras satisfying $(\io \ot \mu)\deh_2 = \te \mu$.
\item There exists an $\cR \in H^1(\matched,\T)$ such that
$$I = \cR (\cL(G_2) \ot 1) \cR^* \;  ,$$
where we consider $\cR : G_2 \times G_1 \recht \T$ as a unitary
multiplication operator on $\rL^2(G_2 \times G_1)$.
\end{enumerate}
Moreover, the correspondence $\cR \mapsto \cR (\cL(G_2) \ot 1) \cR^*$
defines a bijection between $H^1(\matched,\T)$ and the set of full
coideals of the bicrossed product $(M,\de)$.
\end{proposition}
\begin{proof}
Recall that $H^1(\matched,\T)$ consists of the functions $\cR \in
\rL(G_2 \times G_1,\T)$ such that
$$\bcR(p_2(sg),h) \; \cR(s,gh) \; \bcR(s,g) = 1 \qquad\text{and}\qquad
\cR(t,g) \; \bcR(st,g) \; \cR(s,p_1(tg)) = 1$$
for almost all $(s,t,g,h) \in G_2^2 \times G_1^2$.

Let $\cR$ be such a function and consider $\cR$ as a unitary
(multiplication) operator on $\rL^2(G_2 \times G_1)$. Define
\begin{equation}\label{eq.coideal}
I = \cR (\cL(G_2) \ot 1)\cR^* \; .
\end{equation}
Observe that $I$ is a von Neumann algebra generated by $\cR (\lambda_s
\ot 1) \cR^*$ for $s \in G_2$. Using the second equation satisfied by
$\cR$, one gets
$$\cR(\lambda_s \ot 1) \cR^* = (\lambda_s \ot 1) \pi_1(\cR(s,\cdot))
\; .$$
Hence $I$ is a von Neumann subalgebra of $M$. Using the first relation
satisfied by $\cR$, one verifies immediately that $I$ is a coideal of
$(M,\de)$. We also observe
that $\mu : \cL(G_2) \recht I : \mu(a) = \cR(a \ot 1) \cR^*$ defines
an isomorphism of von Neumann algebras satisfying $(\io \ot \mu)\deh_2
= \te \mu$. Hence, $E_s = \mu(\C \, \lambda_s)$ and $I$ is full.
So, we have proven that there is a map from $H^1(\matched,\T)$ to the
set of coideals satisfying all three conditions in the statement of
the proposition.

We prove that this map is injective. For this, it is
sufficient to prove that, whenever $\cR \in H^1(\matched,\T)$ and \linebreak
$\cR(\cL(G_2) \ot 1) \cR^* = \cL(G_2) \ot 1$, then $\cR=1$. In that case, $(\lambda_s \ot 1)
\pi_1(\cR(s,\cdot)) \in \cL(G_2) \ot 1$ for all $s \in G_2$. Hence,
$\pi_1(\cR(s,\cdot)) \in \C$ and we find a function $\cV : G_2 \recht
\T$ such that $\cR(s,g) = \cV(s)$ for almost all $(s,g) \in G_2 \times
G_1$. Considering $\cR \in \rL(\X_{11},\T)$ and $\cV \in
\rL(\X_{01},\T)$, this means that $\cR = \coborv_1 \cV$. Then,
shifting to additive notation
$$0 = \coborv \cR = (\coborv_0 - \coborv_1 + \coborv_2) \coborv_1 \cV
= \coborv_0 \coborv_1 \cV \; .$$
This implies that $\cV$ is trivial and hence $\cR = 0$ in
$H^1(\matched,\T)$.

Suppose next that $I$ is a full coideal in the bicrossed product
$(M,\de)$. We shall produce an element in $H^1(\matched,\T)$ such that
$I$ is given by Equation \eqref{eq.coideal}. Hence, $I$ will satisfy also
the second statement of the proposition.

We claim that there
exists a Borel measurable map $Y : G_2 \recht \cU(I)$ such that
$\te(Y(s)) = \lambda_s \ot Y(s)$. Here, $\cU(I)$ denotes the unitary
group of $I$ equiped with the strong topology. As such, $\cU(I)$ is a
Polish group. Let $s \in G_2$ and $x \in E_s$. Then, $x^* x, x x^* \in
E_e$. Further, $\C \, 1 \subset E_e$ implying that $E_e = \C \, 1$. This
means that $E_s$ is of the form $\C \, u$, with $u$ a unitary in
$I$. Define the Polish group $K \subset G_2 \times \cU(I)$ as
$$K = \{(s,u) \in G_2 \times \cU(I) \mid \te(u) = \lambda_s \ot u \} \; .$$
Define $\pi : K \recht G_2 : \pi(s,u) = s$. Then, $\pi$ is continuous
and surjective. It follows that there exists a Borel section $G_2
\recht K$ and this proves our claim.

Denote by $M^\te$ the subalgebra of $M$ consisting of $x \in M$
satisfying $\te(x) = 1 \ot x$. Observe that, since $\te$ is the dual
coaction on the crossed product $G_2 \ltimes \rL^\infty(G_1)$, we get
that $M^\te = \pi_1(\rL^\infty(G_1))$. For any $s \in G_2$,
$(\lambda_s^* \ot 1) Y(s) \in M^\te$. Hence, we can take $\cR : G_2
\times G_1 \recht \T$ Borel measurable such that $Y(s) = (\lambda_s
\ot 1) \pi_1(\cR(s,\cdot))$ for all $s \in G_2$.

Since $I$ is a von Neumann algebra, we have $Y(st)^* Y(s) Y(t) \in I
\cap M^\te = \C$. It follows that there is a measurable function $\cV
: G_2 \times G_1 \recht \T$ such that
$$\overline{\cR(st,g)} \; \cR(s,p_1(tg)) \; \cR(t,g) = \cV(s,t)$$
almost everywhere. Considering $\cR \in \rL(\X_{11},\T)$ and $\cV \in
\rL(\X_{02},\T)$, this equation can be rewritten as $\coborh \cR =
\coborv_1 \cV$.

On the other hand, $\de(I) \subset M \ot I$ and
$$\de(Y(s)) = \bigl(\bigl( (\lambda_s \ot 1 \ot 1) (\pi_1 \ot
\io)(X_s)\bigr) \ot 1\bigr) \; (\pi_1 \ot \pi_1)
\de_1\bigl(\cR(s,\cdot)\bigr) \; .$$
Hence, we find that
$$(\lambda_{p_2(sg)} \ot 1) \pi_1\bigl( \cR(s,g \; \cdot \;) \bigr)
\in I$$
for almost all $s,g$. But also $Y(p_2(sg)) = (\lambda_{p_2(sg)} \ot 1)
\pi_1\bigl(\cR(p_2(sg), \cdot\bigr) \in I$. It follows that we can
take a measurable function $\cU : G_2 \times G_1 \recht \T$ such that
$$\cR(p_2(sg),h) \; \overline{\cR}(s,gh) \; \cR(s,g) = \cU(s,g)$$
almost everywhere. Considering $\cR,\cU \in \rL(\X_{11},\T)$, this
equation can be rewritten as $\coborv \cR = \coborv_2 \cU$.

We shift back to additive notation. From the equation $\coborv \cR = \coborv_2 \cU$, it follows that
$(\coborv_0 - \coborv_1 + \coborv_2 - \coborv_3)\coborv_2 \cU =
\coborv \coborv_2 \cU = 0$, which yields that $\coborv_3 (\coborv_0 -
\coborv_1) \cU = 0$. Hence, $(\coborv_0 - \coborv_1) \cU = 0$ and we
can take $\cW \in \rL(\X_{01},\T)$ such that $\cU = \coborv_0 \cW$.
We then find that
$$\coborv_2 \coborv_0 \coborh \cW = \coborh \coborv_2 \cU = \coborh
\coborv \cR = \coborv \coborh \cR = \coborv \coborv_1 \cV
= (\coborv_0 - \coborv_1 + \coborv_2) \coborv_1 \cV = \coborv_0
\coborv_1 \cV = \coborv_2 \coborv_0 \cV \; .$$
Hence, $\cV = \coborh \cW$. Define $\widetilde{\cR} = \cR - \coborv_1
\cW$. Then, $\coborh \cRtil = 0$ as well as $\coborv \cRtil = 0$. So,
$\cRtil \in H^1(\matched,\T)$. Moreover, the coideal defined by $\cRtil$ is
generated by the operators
$$(\lambda_s \ot 1) \pi_1(\cRtil(s,\cdot)) = \overline{\cW}(s) \;
(\lambda_s \ot 1) \pi_1(\cR(s,\cdot)) = \overline{\cW}(s) Y(s) \in I \; .$$
Writing $Z(s) = \cRtil (\lambda_s \ot 1) \cRtil^* \in I$, we get a
strongly continuous homomorphism $s \mapsto Z(s)$ satisfying
$\te(Z(s)) = \lambda_s \ot Z(s)$. From Landstad's theorem (see
\cite{Landstad}), it follows that $I$ is a crossed product von Neumann
algebra and in particular, that $I$ is generated by $I^\te$ and
$\{Z(s) \mid s \in G_2\}$. Now, $I^\te = E_e = \C$ and we conclude
that the coideal defined by $\cRtil$ is precisely $I$.
\end{proof}

\begin{remark}
In the previous proposition we described the full coideals of
$(M,\de)$ in terms of $H^1(\matched,\T)$.  There is no hope to
describe \emph{all} coideals of $(M,\de)$ using cohomological
data. This would come down to describe, in the classical setting, all
closed subgroups of $G$ in terms of $G_1$ and $G_2$ if
$$e \longrightarrow G_1 \longrightarrow G \longrightarrow G_2
\longrightarrow e$$
is an exact sequence. It is clear that we can only describe in such a
way a closed subgroup of $G$ whose image in $G_2$ is closed.

In the previous proposition, we characterized in a quantum setting the
closed subgroups $H$ of $G$ such that the homomorphism from $H$ to
$G_2$ is in fact a bijective homeomorphism. More generally, it is
possible to describe the coideals of the bicrossed product $(M,\de)$ whose image in
$\cL(G_2)$ is, in a certain sense, closed.
\end{remark}

\section{Computational methods} \label{sec.computational}

Let $G_1,G_2 \subset G$ be a matched pair of locally compact
groups. We want to calculate the group of extensions
$H^2(\matched,\T)$. Taking into account the Kac exact sequence
(Corollary \ref{cor.Kacexact}), we have to calculate $H^n(\Gamma,\T)$, for
$n=2,3$ and $\Gamma=G,G_1,G_2$. David Wigner developed in \cite{DW} the
necessary tools.

Let $G$ be a locally compact group. The short exact sequence $0 \recht
\Z \recht \R \recht \T \recht 0$ yields the long exact cohomology
sequence
$$\cdots \recht H^n(G,\Z) \recht H^n(G,\R) \recht H^n(G,\T) \recht
H^{n+1}(G,\Z) \recht \cdots$$
So, we have to calculate $H^n(G,\R)$ and $H^n(G,\Z)$.

In order to do so, we explain the approach of David Wigner. We
introduce a special class of Polish $\Z$-modules having {\it property F}, which means that for
  every short exact sequence $0 \recht A \recht B \overset{\pi}{\recht}
  C \recht 0$, $\pi$ has the homotopy lifting property for finite
  dimensional paracompact spaces.

We consider $H^n(G,A)$ for
\begin{itemize}
\item locally compact, $\si$-compact groups $G$ of finite dimension (in
  the sense of Lebesgue covering dimension),
\item Polish $G$-modules $A$ with {\it property F}.
\end{itemize}
In such a situation, we define on the locally compact space $G^n$ the
sheaf $\cA^n$ of continuous functions to $A$, i.e.\ $\cA^n(U) = \{ f :
U \recht A \mid f \;\;\text{continuous} \}$ for all $U \subset G^n$
open. Since the usual face operators $\bor_i : G^n \recht
G^{n-1}$, $i=0,\ldots,n$ are continuous, we get a bicomplex
$(C^k(G^n,\cA^n))_{k,n \geq 0}$:
\begin{equation} \label{eq.sheafbicomplex}
\begin{CD}
\vdots @. \vdots @. \\
@A{\coborv}AA @A{\coborv}AA @. \\
C^1(G^0,\cA^0) @>{\coborh}>> C^1(G^1,\cA^1) @>{\coborh}>> \cdots\\
@A{\coborv}AA @A{\coborv}AA @. \\
C^0(G^0,\cA^0) @>{\coborh}>> C^0(G^1,\cA^1)\ @>{\coborh}>> \cdots
\end{CD}
\end{equation}
Here, $C^k(G^n,\cA^n)$ denotes the
semisimplicial $k$-cochains of sheaf cohomology (see \cite{Bredon},
Chapter II, end of Section 2), $\coborv$ is the coboundary of sheaf
cohomology and $\coborh$ is the coboundary of group cohomology. {\it Wigner
shows that the cohomology of the total bicomplex is precisely
$H(G,A)$.}

One can write an explicit isomorphism as follows: define
$\cA\borel^n$ to be the sheaf on $G^n$ defined by $\cA^n(U) = \{ f :
U \recht A \mid f \;\;\text{Borel measurable} \}$ for all $U \subset G^n$
open. We get a bicomplex $(C^k(G^n,\cA\borel^n))$. Using the inclusion
$\cA^n \hookrightarrow \cA\borel^n$ we get a
cochain transformation $(C^k(G^n,\cA^n)) \recht
(C^k(G^n,\cA\borel^n))$ which is a total cohomology isomorphism. Since the sheafs
$\cA\borel^n$ are flabby, the total cohomology of the second bicomplex
only lives in the first
row and gives precisely the Borel cohomology of
$G$, defined using the complex $(\funborel(G^n,A))_n$.

Suppose now that $A$ is contractible. Then, the
sheafs $\cA^n$ are flabby and the cohomology of the total bicomplex \eqref{eq.sheafbicomplex}
only lives in the first row. By definition, we get that $H(G,A) =
H\cont(G,A)$. Suppose moreover that $A$ is a finite-dimensional vector
space and that $G$ is a Lie group with finitely many connected
components. Let $K \subset G$ be a maximal compact subgroup. Using
results of Hochschild \& Mostow and Van Est, we get that $H(G,A) =
H\cont(G,A) = H(\lieg,K,A)$, the Lie algebra cohomology relative to
$K$ as defined in \cite{Guich}, Chapitre II, n$^\text{o}$ 3.6.

In our examples, we will look at $A=\R$ and $G$ will be
low-dimensional, so that $H(\lieg,K,A)$ is perfectly computable.

On the opposite suppose that $A$ is a discrete
$G$-module. Let $EG \recht BG$ be a principal universal $G$-bundle with
paracompact base (see \cite{husemoller}). Let $\cB^n$ denote the sheaf of continuous (hence,
locally constant) functions
on $EG \times G^n$, i.e.\ $\cB^n(U) = \{ f :
U \recht A \mid f \;\;\text{continuous} \}$ for all $U \subset EG
\times G^n$ open. Define the face operators $\borh_i : EG \times
G^n \recht EG \times G^{n-1}$, $i=0,\ldots,n$ by the formula
$$\borh_i(x,g_1,\ldots,g_n) =
\begin{cases}
(x\cdot g_1,g_2,\ldots,g_n) \quad \text{if}\;\; i=0 , \\
(x,g_1,\ldots,g_i g_{i+1},\ldots,g_n) \quad\text{if}\;\;
i=1,\ldots,n-1, \\
(x,g_1,\ldots,g_{n-1}) \quad\text{if}\;\; i=n \; .
\end{cases}
$$
Proceeding as above, we get a bicomplex $(C^k(EG \times
G^n,\cB^n))_{k,n \geq 0}$. Using the projection $EG \times G^n \recht
G^n$ we get a cochain transformation $(C^k(G^n,\cA^n)) \recht (C^k(EG \times
G^n,\cB^n))$. Since $EG$ is contractible and $A$ discrete, the
homotopy invariance of sheaf cohomology implies that this cochain
transformation induces a cohomology isomorphism in every
column. Hence, we get an isomorphism in total cohomology. Next (and
this argument is slightly delicate, see page 92 in Wigner's paper
\cite{DW}), we can use local sections for the
map $EG \recht BG$ to show that the total cohomology of the bicomplex
$(C^k(EG \times G^n,\cB^n))$ only lives in the first column. As a
conclusion, we get that $H(G,A)$ is precisely the sheaf cohomology of
$BG$ with coefficients in the locally constant sheaf $A$ on $BG$.

As an application, we have the following result.

\begin{proposition} \label{prop.mappingcone}
Let $G_1,G_2 \subset G$ be a matched pair and suppose that $G$ is
finite dimensional. Let $EG \recht BG$ be a principal universal
$G$-bundle with paracompact base. We get principal universal
$G_i$-bundles $EG \recht BG_i$
and continuous maps $BG_i \recht BG$.
The Kac cohomology
$H(\matched, A)$ with coefficients in a discrete, trivial $G$-module
$A$ is the singular cohomology of the mapping cone of $BG_1 \sqcup BG_2
\recht BG$ with coefficients in $A$.
\end{proposition}

\begin{proof}
Denote, on all kinds of spaces, the constant sheaf $A$ by $\cA$ and
the sheaf of Borel functions to $A$ by $\cA\borel$. Consider the
diagram in Equation \eqref{eq.diagram} below.
\begin{figure}[ht]
\begin{equation}\label{eq.diagram}
\begin{CD}
(\rL(\X_{pq},A))_{p+q=n} @>>> (\rL(G_1^n,A) \oplus \rL(G_2^n,A))_n \\
@A{(1)}AA @AA{(1)}A \\
(\funborel(G^n,A))_n @>>> (\funborel(G_1^n,A) \oplus
\funborel(G_2^n,A))_n \\
@V{(2)}VV @VV{(2)}V \\
(C^p(G^q,\cA\borel))_{p+q=n} @>>> (C^p(G_1^q,\cA\borel) \oplus
C^p(G_2^q,\cA\borel))_{p+q=n} \\
@A{(2)}AA @AA{(2)}A \\
(C^p(G^q,\cA))_{p+q=n} @>>> (C^p(G_1^q,\cA) \oplus
C^p(G_2^q,\cA))_{p+q=n} \\
@V{(3)}VV @VV{(3)}V \\
(C^p(EG \times G^q,\cA))_{p+q=n} @>>> (C^p(EG \times G_1^q,\cA) \oplus
C^p(EG \times G_2^q,\cA))_{p+q=n} \\
@A{(4)}AA @AA{(4)}A \\
(C^n(BG, \cA))_n @>>> (C^n(BG_1 \sqcup BG_2,\cA))_n
\end{CD}
\end{equation}
\caption{Diagram for the proof of Proposition \ref{prop.mappingcone}}
\end{figure}
The arrows labeled with $(1)$ are defined by Equation
\eqref{eq.isomorphism}. The arrows labeled with $(2)$ are the obvious
ones. Both are cohomology isomorphisms by the Buchsbaum's criterion
(see Preliminaries).
The arrows labeled with $(3)$ are
cohomology isomorphisms by homotopy invariance of sheaf cohomology, as
explained above. Finally, the arrows labeled with $(4)$ are cohomology
isomorphisms by Wigner's argument on page 92 of \cite{DW}.

The cohomology of the mapping cone of the cochain transformation on
the first line of the above diagram is by definition the Kac
cohomology $H(\matched,A)$. The cohomology of the mapping cone of the cochain
transformation on the last line of the above diagram is precisely the
singular cohomology of the mapping cone of $BG_1 \sqcup BG_2 \recht BG$ with coefficients
in $A$. So, both cohomologies are isomorphic.
\end{proof}

\section{Examples}

\subsection{Matched pairs of low-dimensional Lie groups}

In \cite{VV2}, matched pairs of real Lie algebras $\lieg_1,\lieg_2
\subset \lieg$ have been classified for $\dim \lieg_1 \leq 2$ and
$\lieg_2 = \R$. All these matched pairs have exponentiations to
matched pairs of Lie groups, in the sense of Definition
\ref{def.matched}. The exponentiations that are as connected as
possible (but not always connected) have been given explicitly in
\cite{VV2}. The group of extensions was calculated on
the Lie algebra level and explicit exponentiations to cocycles
satisfying Equation \eqref{eq.cocycles} have been given. Nevertheless,
there were not the necessary tools at hand to prove that these
cocycles really represented exactly the elements of the group of extensions in the
sense of Definition \ref{def.Kaccohomology}.

Using the methods developed above, we will compute correctly the
groups of extensions for the matched pairs of \cite{VV2}, Theorem
5.1.

Remark that if $G_1,G_2$ and $G$ are connected and have no
non-trivial compact subgroups, all cohomologies with coefficients in
$\Z$ vanish and those with coefficients in $\R$ reduce to Lie algebra
cohomology. So, the group of extensions $H^2(\matched,\T)$ is
isomorphic to the group of Lie bialgebra extensions for the
matched pair $\lieg_1,\lieg_2 \subset \lieg$, as defined in
\cite{Mas}. In particular, the group of extensions for
the cases 1~--~3 of \cite{VV2}, Theorem 5.1 are either $\R$ or $0$, as
described in \cite{VV2}, Proposition 6.2.

Next, we have to look at case 4 of \cite{VV2}, Theorem 5.1.
A first
class of matched pairs depends on a parameter $d \not\in \{0,1\}$. We have
\begin{align*}
G &= \Bigl\{ \begin{pmatrix} s & 0 & x \\ 0 & s^d & y \\ 0 & 0 & 1
\end{pmatrix} \Big| \; s \neq 0, x,y \in \R \Bigr\}
\quad\text{where}\quad s^d = \operatorname{Sgn}(s) \, |s|^d \; , \\
G_1 &= \Bigl\{ \begin{pmatrix} s & 0 & 0 \\ 0 & s^d & y \\ 0 & 0 & 1
\end{pmatrix} \Big| \; s \neq 0, y \in \R \Bigr\} \; , \quad\quad
G_2 = \Bigl\{ \begin{pmatrix} s & 0 & s-1 \\ 0 & s^d & 0 \\ 0 & 0 & 1
\end{pmatrix} \Big| \; s \neq 0 \Bigr\} \; .
\end{align*}
Using the results of
Section \ref{sec.computational}, we get that $H^2(G_1,\T) =
H^2(G_2,\T) = 0$, that $H^3(G_1,\R) = H^3(G_2,\R) = 0$ and that the
arrow $H^4(G,\Z) \recht H^4(G_1,\Z)$ is an isomorphism. Further, we
have
$$H^3(G,\R) = \begin{cases} \R \quad\text{if}\;\; d=-1 \; , \\
0 \quad\text{otherwise} \; .\end{cases}$$
All this information, together with the Kac exact sequence, gives that
the group of extensions is $\R$ for $d=-1$ and $0$ otherwise. An
immediate verification gives that, for the case $d=-1$, the cocycles
are exactly given by Proposition 6.3 in \cite{VV2}.

A second class of matched pairs depends on a parameter $b \in \R$. We
have
\begin{align*}
G &= \Bigl\{ \begin{pmatrix} s & b s \log |s| & x \\ 0 & s & y \\ 0 & 0 & 1
\end{pmatrix} \Big| \; s \neq 0, x,y \in \R \Bigr\}
\; , \\
G_1 &= \Bigl\{ \begin{pmatrix} s & b s \log |s| & x \\ 0 & s & 0 \\ 0 & 0 & 1
\end{pmatrix} \Big| \; s \neq 0, x \in \R \Bigr\} \; , \quad\quad
G_2 = \Bigl\{ \begin{pmatrix} s & b s \log |s| & 0 \\ 0 & s & s-1 \\ 0 & 0 & 1
\end{pmatrix} \Big| \; s \neq 0 \Bigr\} \; .
\end{align*}
With exactly the same reasoning as above, we find that the group of
extensions is trivial, since now $H^3(G,\R) =0$.

Finally, consider three more interesting cases.
Before writing the corresponding Lie groups $G_1,G_2 \subset G$, we
compute and write generators for $H^4(\T,\Z)$ and $H^4(\Z/2\Z \ltimes
\T,\Z)$, where $\Z / 2\Z$ acts by taking the inverse. It is well known
that $B\T = P\C^\infty$, the infinite dimensional projective
plane. Under the cup product $H(P\C^\infty,\Z) = \Z[X]$, the
polynomial ring over $\Z$. In particular, if $x$ is a generator for
$H^2(P\C^\infty,\Z)$, $x \cup x$ will be a generator for $H^4(P\C^\infty,\Z)$.
Also the measurable cohomology $H(G,\Z)$ has a cup product: if $\al
\in \rL(G^n,\Z)$ and $\be \in \rL(G^m,\Z)$ are cocycles, we define
$$(\al \cup \be)(g_1,\ldots,g_{n+m}) = \al(g_1,\ldots,g_n)
\be(g_{n+1},\ldots,g_{n+m}) \; .$$
One can verify that the isomorphism $H(G,\Z) \cong H(BG,\Z)$ preserves
cup products. We conclude that, if $\om$ is a generator for $H^2(G,\Z)
\cong \Z$, then $\om \cup \om$ is a generator for $H^4(G,\Z)$. But,
a generator $\om$ for $H^2(\T,\Z)$ is well known, since it comes from the
extension $0 \recht \Z \recht \R \recht \T \recht 0$. We make the
following choice that will be useful later:
$$\text{for} \;  -\pi \leq s,t < \pi \; , \; \text{we define}\quad \om(\exp(it),\exp(is)) = \begin{cases} 1 \quad\text{if}\;\; s+t \geq \pi
  \; , \\ 0 \quad\text{if}\;\; - \pi \leq s+t < \pi \; , \\
-1 \quad\text{if}\;\; s+t < - \pi \; .
\end{cases}$$
An explicit generator for $H^4(\T,\Z)$ is then given by
\begin{equation}\label{eq.alpha}
\al(\lambda_1,
\lambda_2,\lambda_3,\lambda_4)=\om(\lambda_1,
\lambda_2) \, \om (\lambda_3,\lambda_4) \; .
\end{equation}

We turn next to $\Z/2\Z \ltimes \T$. Consider more generally a
semi-direct product $G:=\Gamma \ltimes K$, with $\Gamma$ discrete. Using
the subgroups
$\Gamma$ and $K$, we have a matched pair and hence, we have the
bicomplex in Equation \eqref{eq.bigbicomplex} to compute the
measurable cohomology $H(G,A)$. Associated with this
bicomplex is a spectral sequence (the Lyndon-Hochschild-Serre spectral sequence), which
makes perfectly sense because $\Gamma$ is discrete. We have $E^2_{p,q}
= H^p(\Gamma,H^q(K,A))$, where we consider $H^q(K,A)$ as a discrete
$\Gamma$-module. In our concrete case, with $\Gamma = \Z/2\Z$ and $K =
\T$, we find $E_2^{4,0}=\Z/2\Z$, $E_2^{3,1} = E_2^{2,2} = E_2^{1,3} =
0$ and $E_2^{0,4} = \Z$. Since the action of $\Z/2\Z$ is trivial on
the element $\al \in \rL(\T^4,\Z)$,
we can conclude that $H^4(\Z/2\Z \ltimes \T,\Z) \cong
\Z/2\Z \oplus \Z$ and more specifically, the restriction homomorphism $H^4(\Z/2\Z \ltimes \T,\Z)
\recht H^4(\Z/2\Z,\Z) \oplus H^4(\T,\Z)$ is an isomorphism.

Consider the following matched pair.
\begin{align*}
G &= \{X \in M_2(\R) \mid \det X = \pm 1 \} \mod \{\pm 1\} \; ,
\\ G_1 & = \Bigl\{ \begin{pmatrix} |a| & 0 \\ x & \frac{1}{a} \end{pmatrix}
\Big| \; a \neq 0 , x \in \R \Bigr\} \mod \{\pm 1\} \; , \quad
G_2 = \Bigl\{ \begin{pmatrix} |s| & \frac{1}{2}(|s| - \frac{1}{s}) \\
  0 & \frac{1}{s} \end{pmatrix} \Big| \; s \neq 0 \Bigr\} \mod \{\pm
1\} \; .
\end{align*}

We get that $H^3(G,\R) = 0$ and $H^2(G_1,\T)=H^2(G_2,\T) = 0$. We can
then easily conclude from the Kac exact sequence that the sequence
$0 \recht H^2(\matched,\T) \recht H^4(G,\Z) \recht H^4(G_1,\Z) \oplus
H^4(G_2,\Z)$ is exact. Consider in the obvious way $\Z/2\Z \subset
G_1$ and $\Z/2\Z \subset G_2$. Both embeddings of $\Z/2\Z$ in $G$ are conjugate in $G$. Since conjugation by an element of $G$ acts trivially on
$H^4(G,\Z)$, we conclude that the sequence $0 \recht H^2(\matched,\T)
\recht H^4(G,\Z) \recht H^4(G_1,\Z)$ is exact. Using the maximal
compact subgroup $\Z/2\Z \ltimes \T \subset G$, it follows from this
that $0 \recht H^2(\matched,\T)
\recht H^4(\Z/2\Z \ltimes \T,\Z) \recht H^4(\Z/2\Z,\Z)$ is exact. From
the remarks above, we conclude that $H^2(\matched,\T) \recht
H^4(\T,\Z)$ is an isomorphism. Hence, the group of extensions is $H^2(\matched,\T)
=\Z$. In \cite{VV2}, Proposition 6.3, explicit pairs of $2$-cocycles
$(\cU_n,1)_{n \in \Z}$ have been constructed for the matched pair
$G_1,G_2 \subset G$. Using Remark \ref{rem.inverse} and the final
formula on page 171 as well as the first formula on page 172 of
\cite{VV2}, we can check that the image of $(\cU_n,1)$ under the
isomorphism $H^2(\matched,\T) \recht
H^4(\T,\Z)$ is precisely $-4 n \al$, where $\al$ is defined by Equation
\eqref{eq.alpha}. Once one is able to perform such a calculation, one
can see as well how to change the function $f$ in \cite{VV2} in order
to get exactly the whole of $H^2(\matched,\T)$.

Next, consider the matched pair with
$$G=\operatorname{PSL}_2(\R) \; , \quad
G_1 = \Bigl\{ \begin{pmatrix} a & x \\ 0 & \frac{1}{a} \end{pmatrix}
\Big| \; a > 0, x \in \R \Bigr\} \mod \{\pm 1\} \; , \quad
G_2 = \Bigl\{ \begin{pmatrix} 1 & 0 \\ s & 1 \end{pmatrix}
\Big| \; s \in \R \Bigr\} \mod \{\pm 1\} \; .$$
Considering the maximal compact subgroup $\T \subset G$,
we get immediately that $H^2(\matched,\T) \recht H^4(\T,\Z)$ is an
isomorphism. The explicit cocycles $(\cU_n,1)$ found in \cite{VV2} are
mapped to $2 n \al$ under this isomorphism. Again, it is not hard to
find then explicit formulas for cocycles giving exactly the whole of
$H^2(\matched,\T)$.

Finally, we look at the matched pair
\begin{align*}
G=\operatorname{PSL}_2(\R) \; , \quad G_1 & =
\Bigl\{ \begin{pmatrix} a & x \\ 0 & \frac{1}{a} \end{pmatrix}
\Big| \; a > 0, x \in \R \Bigr\} \mod \{\pm 1\} \; , \\ G_2 &=
\Bigl\{ \begin{pmatrix} \cos t & \sin t \\ - \sin t & \cos t \end{pmatrix}
\Big| \; t \in \R \Bigr\} \mod \{\pm 1\} \; .
\end{align*}
From the Kac exact sequence, it
follows that the sequence $0 \recht H^2(\matched,\T) \recht H^4(G,\Z) \recht
H^4(G_2,\Z)$ is exact. Since $H^4(G,\Z) \recht
H^4(G_2,\Z)$ is an isomorphism, we get $H^2(\matched,\T)=0$. This
explains why it is impossible in \cite{VV2} to exponentiate the
cocycles from the Lie algebra to the Lie group level.

In the final matched pair above, we can take the associated matched
pair of Lie algebras $\lieg_1,\lieg_2 \subset \lieg$. If we take
$\widetilde{G}$ to be the connected, simply connected Lie
group with Lie algebra $\lieg$ and $\widetilde{G}_1,\widetilde{G}_2$
to be the connected Lie subgroups with Lie algebras $\lieg_1,\lieg_2$,
we get a matched pair such that
$\widetilde{G},\widetilde{G}_1,\widetilde{G}_2$
are connected and without compact
subgroups. So, we conclude as above that
$H^2(\matched,\T)$ is isomorphic to the group of Lie bialgebra
extensions, i.e.\ $\R$ as stated in Proposition 6.2 of \cite{VV2}.

\subsection{Some other examples}

Let $G$ be a semi-simple Lie group with finite center, such that in
its Iwasawa decomposition $G=KAN$, $K$ is a maximal compact
subgroup. Writing $G_1=K$ and $G_2 = AN$, we get a matched pair of Lie
groups. Since $G_2$ is contractible, we get that $H^n(G_2,\Z)=0$ and
$H^n(G,\Z) \recht H^n(G_1,\Z)$ is an isomorphism for all $n \geq
1$. It follows from the Kac exact sequence that $H^n(\matched,\Z) = 0$
for $n \geq 1$. Hence, $H^n(\matched,\R) \recht H^n(\matched,\T)$ is
an isomorphism. But,
$$H^2(G,\R) \recht H^2(G_1,\R) \oplus H^2(G_2,\R) \recht H^2(\matched,\R)
\recht H^3(G,\R) \recht H^3(G_1,\R) \oplus H^3(G_2,\R)$$
is an exact sequence of vector spaces. Moreover, $H^n(G_1,\R) = 0$ for
$n \geq 1$. Writing $\lies = \liea \oplus \lien$ for the Lie algebra
of $G_2$ and $\liek$ for the Lie algebra of $K$,
it
follows that the group of extensions is isomorphic with
$$H^2(\matched,\T) \cong \operatorname{Coker}\bigl(
H^2(\lieg,\liek,\R) \recht H^2(\lies,\R) \bigr) \oplus
\operatorname{Ker}\bigl(H^3(\lieg,\liek,\R) \recht H^3(\lies,\R)
\bigr) \; .$$
If we take $G= \operatorname{SL}_2(\R)$, we immediately get that the
group of extensions is trivial. For $G= \operatorname{SL}_2(\C)$, a
direct calculation gives $H^2(\lies,\R) = H^3(\lies,\R) = 0$, while
$H^3(\lieg,\liek,\R) = \R$. So, in this case $H^2(\matched,\T) =
\R$. Observe that we deal precisely with the matched pair considered
by Majid in \cite{Maj1}.

In \cite{BSV}, a class of matched pairs is constructed as follows: let
$\cA$ be a locally compact ring such that $\cA \setminus \cA^\star$
has (additive) Haar measure zero, where $\cA^\star$ is the group of
units. Let $G$ be the group of affine transformations of $\cA$, i.e.\
the locally compact group with underlying space $\cA^\star \times \cA$
and product $(a,x)\cdot(b,y) = (ab,x+ay)$. The subgroups $G_1,G_2$ will be
the subgroups of transformations fixing $0$ and $-1$ respectively. This
means that $G_1$ and $G_2$ consist of the elements $(a,0)$ and
$(b,b-1)$ for $a,b \in \cA^\star$.

If we take $\cA = \R$, we easily get that $H^2(\matched,\T)=0$. If
$\cA = \C$, we observe that in the commutative diagram
\begin{equation*}
\begin{CD}
H^3(G,\T) @>>> H^3(G_1,\T) \\
@VVV @VVV \\
H^3(\T,\Z) @>>> H^3(\T,\Z)
\end{CD}
\end{equation*}
the vertical arrows and the lower horizontal arrow are
isomorphisms. Since $H^2(G_1,\T) = H^2(G_2,\T) = 0$, we can conclude
that $H^2(\matched,\T)=0$.

\end{document}